\documentclass[11pt,english]{article}
\pdfoutput=1
\usepackage[T1]{fontenc}
\usepackage[latin9]{inputenc}
\usepackage{geometry}
\usepackage{color}
\usepackage{units}
\usepackage{amsmath}
\usepackage{amssymb}
\usepackage{stmaryrd}

\makeatletter

\providecommand{\tabularnewline}{\\}

\usepackage{latexsym}
\usepackage{amsfonts}
\usepackage{amsthm}
\usepackage{mathrsfs}
\@ifundefined{definecolor}
 {\usepackage{color}}{}
\@ifundefined{definecolor}{\usepackage{color}}{}

\newcommand{\A}{\mathbf{A}}

\newcommand{\diag}{\text{diag}}

\newcommand{\U}{\mathbf{U}}
\newcommand{\V}{\mathbf{V}}

\theoremstyle{definition}

\theoremstyle{remark}

\usepackage{babel}
\usepackage{nicefrac}
\usepackage{babel}
\usepackage{sagetex}

\makeatother

\usepackage{babel}
\begin{document}
\title{A symmetrization approach to hypermatrix SVD}
\author{Edinah K. Gnang \thanks{Department of Applied Mathematics and Statistics, Johns Hopkins University,
email: egnang1@jhu.edu}, Fan Tian \thanks{Department of Applied Mathematics and Statistics, Johns Hopkins University,
email: ftian4@jhu.edu}}
\maketitle
\begin{abstract}
We propose a new hypermatrix singular value decomposition based upon
the spectral decomposition of the symmetric products of transposes.
\end{abstract}

\section{Introduction}

One of the most fruitful ideas in matrix theory is that of matrix
decomposition or canonical form. Of the many matrix canonical forms
discussed in the literature, the Singular Value Decomposition (or
SVD for short), is by far the most widely used. Recall that for an
arbitrary $\mathbf{A}\in\mathbb{C}^{n\times n}$, the SVD of $\mathbf{A}$
is expressed by 
\begin{equation}
\begin{array}{c}
\mathbf{A}=\left(\U\sqrt{\diag\boldsymbol{\sigma}\left(\A\right)}\right)\left(\sqrt{\diag\boldsymbol{\sigma}\left(\A\right)}\V\right)\\
\text{such that }\\
\mathbf{U}\mathbf{U}^{*}=\mathbf{I}_{n}=\mathbf{V}^{*}\mathbf{V}
\end{array}\label{SVD constraints}
\end{equation}
Calculating the SVD consists of finding the eigenvalues and eigenvectors
of the Hermitian products of $\mathbf{A}$ and $\mathbf{A}^{*}$.
Important infomation about the matrix $\mathbf{A}$ is obtained through
decomposition such as the matrix rank, the orthornormal basis vectors
and the diagonal matrix of the scaling values, all of which are useful
to be exteneded to higher dimensions. Over the past decades, considerable
progress has been made in generalizing the matrix SVD to higher order
hypermatrices. Two predominat approaches to hypermatrix canonical
forms are now well established as the CANDECOMP-PARAFAC (CP) model
\cite{Carroll1970AnalysisOI,Harshman1970FoundationsOT} and the Tucker
model \cite{Tucker1966}, where the former is a special case to the
later. Based on Tucker model, De Lathauwer, De Moor, and Vandewalle
poineered a multilinear generalization of the matrix SVD to hypermatrices
in \cite{doi:10.1137/S0895479896305696}, namely the Higher-Order
Singular Value Decomposition (HOSVD). The classical models of CP and
Tucker or HOSVD generally express the decompostion of a hypermatrix
as a sum of outer products of vectors, also referred to as the $n$-mode
product in the form of ``hypermatrix times matrices'' \cite{doi:10.1137/07070111X}.
In particular, the $n$-mode product enables the hypermatrix SVD through
performing matrix SVDs following the mode-$n$ flattening (unfolding)
of the original hypermatrix into matrices, and then assemble results
into a hypermatrix of the same order. One of the advantages of the
classical models and the method of HOSVD is that the obtained results
guarantee orthogonality to some extent: the singular vectors are entries
of orthogonal matrices, and the core hypermatrix coordinating singular
values meets a property of all-orthogonality that is a relaxation
to the diagonality property in the matrix SVD. Thorough discussions
on the classical methods and applications have been reviewed in \cite{doi:10.1137/07070111X}.
Other more recent studies also explored alternative representations
of a hypermatrix SVD as a sum of outer products of matrices, which
is a generalization based on a different hypermatrix multiplication
scheme in the form of ``hypermatrix times hypermatrix'' \cite{Kilmer2008ATG,KILMER2011641}. 

While the aforementioned generalizations to higher-order SVD have
been widely used in applications, they often reduce the problems to
matrix SVDs through the folding and unfolding schemes. By contrast
to the matrix case, such higher-order SVD methods do not stem from
a hypermatrix formulation of the spectral theorem. Recent works in
\cite{Gnang2017238,GER} motivated by the generalization of the spectral
theorem to hypermatrices suggest new ways to extend matrix SVD to
hypermatrix SVD while retaining the link to the spectra. In the present
note, we discuss in analogy to matrix SVD the new approach to obtain
orthogonal hypermatrices and diagonal scaling hypermatrix via spectral
decompostion of symmetric products of transposes. Our work is based
on the the Bhattacharya-Mesner algebra (BM algebra) introduced in
\cite{MB94,Gnang2017238,GNANG2020391}, which has enabled the generalization
of many important matrix concepts including the rank, inverse, and
spectral decompostions to hypermatrices. In addition to the hypermatrix
SVD, we also expand the list of concepts to the BM algebra to include
definitions of tensorial orbits and invariants of hypermarices, and
hypermatrix orthorgonality and unitarity.

\section{Overview of the Bhattacharya-Mesner (BM) algebra.}

\emph{Hypermatrices} are multidimensional matrices. More precisely,
a hypermatrix is a finite multiset whose elements (called entries)
are indexed by members of some fixed Cartesian product of the form
\[
\left\{ 0,\cdots,n_{0}-1\right\} \times\left\{ 0,\cdots,n_{1}-1\right\} \times\cdots\times\left\{ 0,\cdots,n_{m-1}-1\right\} .
\]
Such a hypermatrix is of order $m$ and of size $n_{0}\times n_{1}\times\cdots\times n_{m-1}$.
A hypermatrix is cubic of side length $n$ if $n_{0}=n_{1}=\cdots=n_{m-1}=n$. 

Hypermatrix algebras arise from natural generalizations of classical
matrix notions and algorithms \cite{MB94,GKZ,RK,GER,Gnang2017238,MB90}.
The important distinction between hypermatrices and tensors closely
mirrors the distinction between matrices and abstract linear transformations.
Recall that an abstract linear transformation specified over finite
dimensional $\mathbb{K}$-vector spaces is identified with a matrix
orbit. For instance, let $\mathbf{M}\in\mathbb{K}^{m\times n}$ be
associated with some abstract linear transformation specified relative
to the standard basis for $\mathbb{K}^{n\times1}$ and $\mathbb{K}^{1\times m}$.
The tensorial orbit of the linear transformation (accounting for all
possible coordinate changes) is the matrix set 
\[
\left\{ \mathbf{A}\cdot\mathbf{M}\cdot\mathbf{B}\,:\,\begin{array}{c}
\mathbf{A}\in\text{GL}_{m}\left(\mathbb{K}\right)\\
\text{and}\\
\mathbf{B}\in\text{GL}_{n}\left(\mathbb{K}\right)
\end{array}\right\} .
\]
A matrix property common to every member of a tensorial orbit is a
\emph{tensorial invariant}.

Classically, third order hypermatrices in $\mathbb{K}^{m\times n\times p}$
arise from tensorial orbits induced by the action of various appropriate
subgroups of the general linear group on canonical embeddings of $\mathbb{K}$-vector
spaces : $\mathbb{K}^{m\times1\times1}$, $\mathbb{K}^{1\times n\times1}$
and $\mathbb{K}^{1\times1\times p}$ respectively. Incidentally, classical
tensorial invariants such as the rank and singular values are defined
by analogy to their matrix counterparts.

Hypermatrix multiplication, named the \emph{Bhattacharya-Mesner product}
(BM-product), is a generalization to the matrix multiplication \cite{MB90,MB94}.
Occasionally, the product of a conformable matrix pair
\[
\mathbf{A}\in\mathbb{K}^{m\times\textcolor{red}{\ell}},\,\mathbf{B}\in\mathbb{K}^{\textcolor{red}{\ell}\times n},
\]
can be written using the BM-product notation as $\text{Prod}\left(\mathbf{A},\mathbf{B}\right)$
for consistency and such a product is specified entry-wise by 

\[
\text{Prod}\left(\mathbf{A},\mathbf{B}\right)\left[i,j\right]=\sum_{\textcolor{red}{0\le t<\ell}}\mathbf{A}\left[i,\textcolor{red}{t}\right]\,\mathbf{B}\left[\textcolor{red}{t},j\right],\quad\forall\:\begin{cases}
\begin{array}{c}
0\le i<m\\
0\le j<n
\end{array}\end{cases}.
\]

Similarly, the BM-product of a conformable triple of third order hypermatrices
\[
\mathbf{A}\in\mathbb{K}^{m\times\textcolor{red}{\ell}\times p},\,\mathbf{B}\in\mathbb{K}^{m\times n\times\textcolor{red}{\ell}}\mbox{ and }\mathbf{C}\in\mathbb{K}^{\textcolor{red}{\ell}\times n\times p},
\]
is noted $\text{Prod}\left(\mathbf{A},\mathbf{B},\mathbf{C}\right)$
and specified entry-wise by
\[
\text{Prod}\left(\mathbf{A},\,\mathbf{B},\,\mathbf{C}\right)\left[i,j,k\right]=\sum_{\textcolor{red}{0\le t<\ell}}\mathbf{A}\left[i,\textcolor{red}{t},k\right]\,\mathbf{B}\left[i,j,\textcolor{red}{t}\right]\,\mathbf{C}\left[\textcolor{red}{t},j,k\right],\quad\forall\:\begin{cases}
\begin{array}{c}
0\le i<m\\
0\le j<n\\
0\le k<p
\end{array}\end{cases}.
\]

Furthermore, we recall that the \emph{general Bhattacharya-Mesner
product} of a conformable triple 
\[
\mathbf{A}\in\mathbb{K}^{m\times\textcolor{red}{\ell}\times p},\,\mathbf{B}\in\mathbb{K}^{m\times n\times\textcolor{red}{\ell}}\mbox{ and }\mathbf{C}\in\mathbb{K}^{\textcolor{red}{\ell}\times n\times p},
\]
taken with an additional cubic \emph{background hypermatrix} $\mathbf{M}\in\mathbb{K}^{\textcolor{red}{\ell}\times\textcolor{red}{\ell}\times\textcolor{red}{\ell}}$
(similar to metric tensors first introduced in differential geometry
\cite{Ricci1900,gauss1828disquisitiones}) is denoted $\mbox{Prod}_{\mathbf{M}}\left(\mathbf{A},\,\mathbf{B},\,\mathbf{C}\right)\in\mathbb{K}^{m\times n\times p}$
and specified entry-wise by
\begin{equation}
\mbox{Prod}_{\mathbf{M}}\left(\mathbf{A},\,\mathbf{B},\,\mathbf{C}\right)\left[i,j,k\right]=\sum_{\textcolor{red}{0\le t_{0},t_{1},t_{2}<\ell}}\mathbf{A}\left[i,\textcolor{red}{t_{0}},k\right]\mathbf{B}\left[i,j,\textcolor{red}{t_{1}}\right]\mathbf{C}\left[\textcolor{red}{t_{2}},j,k\right]\mathbf{M}\left[\textcolor{red}{t_{0}},\textcolor{red}{t_{1}},\textcolor{red}{t_{2}}\right].\label{General BM product}
\end{equation}

The original BM-product is thus recovered from the general BM-product
by setting the cubic background hypermatrix $\mathbf{M}$ to be equal
to the \emph{Kronecker delta} hypermatrix denoted $\boldsymbol{\Delta}$,
whose entries are specified by 
\[
\boldsymbol{\Delta}\left[i_{0},i_{1},i_{2}\right]=\begin{cases}
\begin{array}{cc}
1 & \mbox{ if }\:0\le i_{0}=i_{1}=i_{2}<n\\
0 & \mbox{otherwise}
\end{array}\end{cases}.
\]
The \emph{general Bhattacharya-Mesner product} of conformable matrices
\[
\mathbf{A}\in\mathbb{K}^{m\times\textcolor{red}{\ell}},\,\mathbf{B}\in\mathbb{K}^{\textcolor{red}{\ell}\times n},
\]
taken with the background matrix $\mathbf{M}\in\mathbb{K}^{\textcolor{red}{\ell}\times\textcolor{red}{\ell}}$
is given by 
\[
\text{Prod}_{\mathbf{M}}\left(\mathbf{A},\mathbf{B}\right)\left[i,j\right]=\sum_{\textcolor{red}{0\le t_{0},t_{1}<\ell}}\mathbf{A}\left[i,\textcolor{red}{t_{0}}\right]\,\mathbf{B}\left[\textcolor{red}{t_{1}},j\right]\mathbf{M}\left[\textcolor{red}{t_{0}},\textcolor{red}{t_{1}}\right],\quad\forall\:\begin{cases}
\begin{array}{c}
0\le i<m\\
0\le j<n
\end{array}\end{cases}.
\]

We further recall that the transpose of an arbitrary hypermatrix $\mathbf{A}\in\mathbb{K}^{m\times n\times p}$,
denoted as $\mathbf{A}^{\top}\in\mathbb{K}^{n\times p\times m}$,
results from a cyclic permutation on the indices and is specified
entry-wise as follows
\[
\mathbf{A}^{\top}\left[i,j,k\right]=\mathbf{A}\left[k,i,j\right].
\]
We adopt the convention
\[
\mathbf{A}^{\top^{2}}:=\left(\mathbf{A}^{\top}\right)^{\top},\ \mathbf{A}^{\top^{3}}:=\left(\mathbf{A}^{\top^{2}}\right)^{\top}=\mathbf{A}.
\]
\[
\implies\mathbf{A}^{\top^{u}}=\mathbf{A}^{\top^{v}}\:\mbox{ if }\:u\equiv v\mod3.
\]
Note that when $\mathbb{K}$ is commutative 
\[
\text{Prod}\left(\mathbf{A},\,\mathbf{B},\,\mathbf{C}\right)^{\top}=\text{Prod}\left(\mathbf{B}^{\top},\,\mathbf{C}^{\top},\,\mathbf{A}^{\top}\right).
\]

\section{Tensorial matrix orbits.}

Let $\mathbb{K}$ denote an arbitrary field (not necessarily commutative)
and let $\text{GL}_{n}\left(\mathbb{K}\right)$ denote the \emph{general
linear group} of invertible $n\times n$ matrices whose entries belong
to $\mathbb{K}$. When investigating matrices, it is of interest to
determine matrix attributes which are independent of the chosen coordinate
system. For this purpose we associate with an arbitrary matrix $\mathbf{M}\in\mathbb{K}^{m\times n}$
a \emph{tensorial orbit} induced by the action on $\mathbf{M}$ of
the group $\text{GL}_{m}\left(\mathbb{K}\right)\times\text{GL}_{n}\left(\mathbb{K}\right)$
as follows

\begin{equation}
\mathcal{T}\left(\mathbf{M}\right):=\left\{ \mathbf{A}\cdot\mathbf{M}\cdot\mathbf{B}\,:\,\begin{array}{c}
\mathbf{A}\in\text{GL}_{m}\left(\mathbb{K}\right)\\
\text{and}\\
\mathbf{B}\in\text{GL}_{n}\left(\mathbb{K}\right)
\end{array}\right\} .\label{Matrix Tensorial Orbit}
\end{equation}
For instance, the tensorial orbit of $\left(\begin{array}{cc}
1 & 1\\
0 & 1
\end{array}\right)$ whose entries are taken from the finite field with two elements denoted
$\mathbb{F}_{2}$ is
\[
\left\{ \left(\begin{array}{cc}
1 & 1\\
0 & 1
\end{array}\right),\,\left(\begin{array}{cc}
1 & 1\\
1 & 0
\end{array}\right),\,\left(\begin{array}{cc}
1 & 0\\
1 & 1
\end{array}\right),\,\left(\begin{array}{cc}
0 & 1\\
1 & 1
\end{array}\right),\,\left(\begin{array}{cc}
1 & 0\\
0 & 1
\end{array}\right),\,\left(\begin{array}{cc}
0 & 1\\
1 & 0
\end{array}\right)\right\} .
\]
In particular the tensorial orbit of a zero matrix is a singleton
\[
\mathcal{T}\left(\mathbf{0}_{m\times n}\right)=\left\{ \mathbf{0}_{m\times n}\right\} .
\]
Recall that 
\[
\forall\:\mathbf{M}\in\text{GL}_{n}\left(\mathbb{K}\right),\quad\mathcal{T}\left(\mathbf{M}\right)=\text{GL}_{n}\left(\mathbb{K}\right).
\]
When $\mathbb{K}=\mathbb{F}_{p^{k}}$ for a prime $p$ we have
\[
\left|\mathcal{T}\left(\mathbf{M}\right)\right|=\prod_{0\le i<n}\left(p^{k\cdot n}-p^{k\cdot i}\right).
\]

The cardinality $\left|\mathcal{T}\left(\mathbf{M}\right)\right|$
is by definition a tensorial invariant, whereas the property of being
symmetric (i.e. $\mathbf{M}=\mathbf{M}^{\top}$) is not in general
a tensorial invariant. Classical matrix attributes well known to be
tensorial invariants include : 
\begin{itemize}
\item The \emph{rank} of $\mathbf{M}\in\mathbb{C}^{m\times n}$ defined
as
\[
\min_{\begin{array}{c}
\mathbf{A}\in\text{GL}_{m}\left(\mathbb{C}\right)\\
\mathbf{B}\in\text{GL}_{n}\left(\mathbb{C}\right)
\end{array}}\left\Vert \mathbf{A}\,\mathbf{M}\,\mathbf{B}\right\Vert _{\ell_{0}}
\]
\item The \emph{nullity} of $\mathbf{M}\in\mathbb{C}^{m\times n}$ defined
as
\[
\text{Dimension of }\left\{ \mathbf{x}\in\mathbb{C}^{1\times m}\;:\;\mathbf{x}\,\mathbf{A}\,\mathbf{y}=0,\:\forall\:\mathbf{y}\in\mathbb{C}^{n\times1}\right\} =\min\left(m,n\right)-\min_{\begin{array}{c}
\mathbf{A}\in\text{GL}_{m}\left(\mathbb{C}\right)\\
\mathbf{B}\in\text{GL}_{n}\left(\mathbb{C}\right)
\end{array}}\left\Vert \mathbf{A}\,\mathbf{M}\,\mathbf{B}\right\Vert _{\ell_{0}}.
\]
\item Singular values of $\mathbf{M}\in\mathbb{C}^{m\times n}$, defined
as multiset of moduli of diagonal entries of any diagonal matrix element
in the sub-orbit of $\mathcal{T}\left(\mathbf{M}\right)$ 
\[
\left\{ \mathbf{A}\,\mathbf{M}\,\mathbf{B}\,:\,\begin{array}{c}
\mathbf{A}\in\text{U}_{m}\left(\mathbb{C}\right)\\
\text{and}\\
\mathbf{B}\in\text{U}_{n}\left(\mathbb{C}\right)
\end{array}\right\} ,
\]
where $\text{U}_{m}\left(\mathbb{C}\right)$ and $\text{U}_{n}\left(\mathbb{C}\right)$
respectively denote the unitary subgroup of GL$_{m}\left(\mathbb{C}\right)$
and GL$_{n}\left(\mathbb{C}\right)$.
\item The eigenvalues of $\mathbf{M}\in\mathbb{C}^{n\times n}$ defined
as
\[
\left\{ \lambda\in\mathbb{C}\::\:0=\det\left(\lambda\mathbf{I}_{n}-\mathbf{A}\,\mathbf{M}\,\mathbf{B}\right)\,:\,\begin{array}{c}
\mathbf{A},\mathbf{B}\in\text{GL}_{n}\left(\mathbb{C}\right)\\
\text{and}\\
\mathbf{I}_{n}=\mathbf{A}\mathbf{B}
\end{array}\right\} .
\]
\end{itemize}

\section{Classical hypermatrix tensorial orbits and their invariants}

Classical hypermatrix tensorial orbits are similar to matrix tensorial
orbits in that they are both resulted from the action of the general
linear group. Hypermatrix tensorial orbits are often simply called
tensors in the literature, which are defined as elements of the tensor
product of vector spaces \cite{silva&lim,lim2013}. The classical
tensorial orbit of the hypermatrix $\mathbf{H}\in\mathbb{K}^{m\times n\times p}$,
resulted from the action of the group $\text{GL}_{m}\left(\mathbb{K}\right)\times\text{GL}_{n}\left(\mathbb{K}\right)\times\text{GL}_{p}\left(\mathbb{K}\right)$
on $\mathbf{H}$ is given by
\[
\mathcal{T}\left(\mathbf{H}\right)\::=\left\{ \mathbf{H}\times_{0}\mathbf{A}\times_{1}\mathbf{B}\times_{2}\mathbf{C}\,:\,\begin{array}{c}
\mathbf{A}\in\text{GL}_{m}\left(\mathbb{K}\right)\\
\mathbf{B}\in\text{GL}_{n}\left(\mathbb{K}\right)\\
\mathbf{C}\in\text{GL}_{p}\left(\mathbb{K}\right)
\end{array}\right\} .
\]
The notation above refers to the $n$-mode product introduced by De
Lathauwer, De Moor, and Vandewalle in \cite{doi:10.1137/S0895479896305696}.
We also note that this notation is equivalent to the multilinear multiplication
in some earlier works denoted as $\left(\mathbf{A},\mathbf{B},\mathbf{C}\right)\cdot\mathbf{H}$
in lieu of $\mathbf{H}\times_{0}\mathbf{A}\times_{1}\mathbf{B}\times_{2}\mathbf{C}$
\cite{silva&lim}.

By analogy to the matrix case, classical third order tensorial invariants
include :
\begin{itemize}
\item The \emph{tensor rank} of $\mathbf{H}\in\mathbb{K}^{m\times n\times p}$
defined as
\[
\min_{\begin{array}{c}
\mathbf{A}\in\text{GL}_{m}\left(\mathbb{K}\right)\\
\mathbf{B}\in\text{GL}_{n}\left(\mathbb{K}\right)\\
\mathbf{C}\in\text{GL}_{p}\left(\mathbb{K}\right)
\end{array}}\left\Vert \mathbf{H}\times_{0}\mathbf{A}\times_{1}\mathbf{B}\times_{2}\mathbf{C}\right\Vert _{\ell_{0}}
\]
\item The \emph{nullity} of $\mathbf{H}\in\mathbb{K}^{m\times n\times p}$
defined as 
\[
\text{Dimension of }\left\{ \left(\mathbf{x},\mathbf{y},\mathbf{z}\right)\in\mathbb{C}^{m}\times\mathbb{C}^{n}\::\:0=\mathbf{H}\times_{0}\mathbf{x}\times_{1}\mathbf{y}\times_{2}\mathbf{z},\,\forall\ \mathbf{z}\in\mathbb{C}^{p}\right\} 
\]
\item Singular values of $\mathbf{H}\in\mathbb{K}^{m\times n\times p}$,
defined as the multiset of moduli super-diagonal entries of diagonal
elements of the tensorial sub-orbit 
\[
\left\{ \mathbf{H}\times_{0}\mathbf{A}\times_{1}\mathbf{B}\times_{2}\mathbf{C}\,:\,\begin{array}{c}
\mathbf{A}\in\text{U}_{m}\left(\mathbb{K}\right)\\
\mathbf{B}\in\text{U}_{n}\left(\mathbb{K}\right)\\
\mathbf{C}\in\text{U}_{p}\left(\mathbb{K}\right)
\end{array}\right\} 
\]
\end{itemize}

\section{Non-classical tensorial orbits and their invariants}

Historically, the study of classical tensorial orbits has been the
predominant approach to investigating hypermatrices \cite{GORDAN1869,Hilbert1890,Ricci1900,gauss1828disquisitiones}.
Unfortunately, two main drawbacks plague the classical tensorial orbits.
The first drawback is conceptual in nature. It results from the fact
that classical tensorial invariants do not suggest a distinct hypermatrix
analog of the general linear group, nor do they suggest any generalization
to hypermatrices of such notions as inverse, nullity, determinant,
spectral decomposition, Rayleigh quotient inequality, resolution of
identity, Parseval identity, unitarity and Fourier transforms. The
second drawback is somewhat related to the first one but is of a computational
nature. Classical tensorial invariant do not suggest any generalization
of classical matrix algorithms such as the rank revealing LU decomposition
and the The Gram-Schmidt orthogonalization process among others. These
drawback have been recently addressed by the proposing new non-classical
tensorial orbits and invariants \cite{GNANG2020391,Gnang2017238}.
For instance, new hypermatrix invariants which extend matrix notions
and algorithms to hypermatrices arise from the BM algebra \cite{GNANG2020391}.
To be more specific, the BM algebra suggests a generalization to higher
order hypermatrices of notions such as inverse and rank so as to enable
the generalization to hypermartrices of the classical Rank Nullity
theorem \cite{GNANG2020391}. On the computational side, the BM approach
also suggest a generalization to hypermatrices of the rank revealing
LU factorization as well as the orthogonalization procedure, and higher
order generalization of the Fourier transforms \cite{GNANG2020391,Gnang2017238}.
The BM algebra also enables a hypermatrix formulation of the spectral
decomposition which we can extend to the symmetrization formulation
of the third order hypermatrix SVD. This latter topic is the main
subject of the present note and will be discussed at length. 

We briefly recall here for the readers' benefit an example of a non-classical
tensorial orbit. Recall that the matrix general linear group over
an arbitrary field $\mathbb{K}$ (possibly non-commutative) is the
matrix set
\begin{equation}
\text{GL}_{m}\left(\mathbb{K}\right)\,:=\left\{ \mathbf{A}\in\mathbb{K}^{m\times m}:\,\exists\:\mathbf{B}\in\mathbb{K}^{m\times m}\text{ s.t. }\mathbf{B}\cdot\mathbf{A}\cdot\mathbf{X}=\mathbf{X},\ \forall\,\mathbf{X}\in\mathbb{K}^{m\times n}\right\} .\label{Definion_of_GLn}
\end{equation}

In contrast to the matrix general linear groups, their third order
hypermatrix analog does not form a group. On the other hand, third
order hypermatrix analog to general linear groups are defined similarly
to Eq. (\ref{Definion_of_GLn}) as follows :
\[
\text{GL}_{m\times n\times p}\left(m\times p\times p,\,p\times n\times p,\,\mathbb{K}\right)\,:=
\]
\begin{equation}
\left\{ \left(\mathbf{A},\mathbf{B}\right)\in\mathbb{K}^{m\times p\times p}\times\mathbb{K}^{p\times n\times p}:\,\exists\:\left(\mathbf{C},\mathbf{D}\right)\in\mathbb{K}^{m\times p\times p}\times\mathbb{K}^{p\times n\times p}\text{ s.t. }\text{Prod}\left(\mathbf{C},\text{Prod}\left(\mathbf{A},\mathbf{X},\mathbf{B}\right),\mathbf{D}\right)=\mathbf{X},\ \forall\,\mathbf{X}\in\mathbb{K}^{m\times n\times p}\right\} .\label{Definition_of_GLmxnxp}
\end{equation}
Just as in the matrix case, third order hypermatrix analog of general
linear groups are defined in terms of \textit{hypermatrix inverse
pairs}. 

Note that over any field (not necessarily commutative) there are subsets
of invertible hypermatrix pairs which do form a group with respect
to the BM product. The simplest example is the third order hypermatrix
analog of the subgroup diagonal matrices. We call such hypermatrices\emph{
scaling hypermatrices.}

A pair $\left(\mathbf{A},\mathbf{B}\right)\in\mathbb{K}^{m\times p\times p}\times\mathbb{K}^{p\times n\times p}$
is an invertible scaling hypermatrix pair if 
\[
\mathbf{A}\left[i,t,k\right]=\begin{cases}
\begin{array}{cc}
\alpha_{it}\in\mathbb{K}\backslash\left\{ 0\right\}  & \mbox{ if }0\le t=k<p\\
0 & \mbox{otherwise}
\end{array} & \forall\,\begin{cases}
\begin{array}{c}
0\le i<m\\
0\le t<p\\
0\le k<p
\end{array}\end{cases}\end{cases},
\]
\[
\mathbf{B}\left[t,j,k\right]=\begin{cases}
\begin{array}{cc}
\beta_{tj}\in\mathbb{K}\backslash\left\{ 0\right\}  & \mbox{ if }0\le t=k<p\\
0 & \mbox{otherwise}
\end{array} & \forall\,\begin{cases}
\begin{array}{c}
0\le t<p\\
0\le j<n\\
0\le k<p
\end{array}\end{cases}\end{cases},
\]
\[
\implies\mbox{Prod}\left(\mathbf{A},\mathbf{X},\mathbf{B}\right)\left[i,j,k\right]=\alpha_{ik}\,\mathbf{X}\left[i,j,k\right]\,\beta_{kj},\quad\forall\:\begin{cases}
\begin{array}{c}
0\le i<m\\
0\le j<n\\
0\le k<p
\end{array}\end{cases}.
\]
The corresponding inverse pair is $\left(\mathbf{C},\mathbf{D}\right)\in\mathbb{K}^{m\times p\times p}\times\mathbb{K}^{p\times n\times p}$
such that 
\[
\mathbf{C}\left[i,t,k\right]=\begin{cases}
\begin{array}{cc}
\alpha_{it}^{-1} & \mbox{ if }0\le t=k\le p\\
0 & \mbox{otherwise}
\end{array} & \forall\,\begin{cases}
\begin{array}{c}
0\le i<m\\
0\le t<p\\
0\le k<p
\end{array}\end{cases}\end{cases},
\]
\[
\mathbf{D}\left[t,j,k\right]=\begin{cases}
\begin{array}{cc}
\beta_{tj}^{-1} & \mbox{ if }0\le t=k\le p\\
0 & \mbox{otherwise}
\end{array} & \forall\,\begin{cases}
\begin{array}{c}
0\le t<p\\
0\le j<n\\
0\le k<p
\end{array}\end{cases}\end{cases}.
\]

Examples of non-classical tensorial orbits associated with $\mathbf{H}\in\mathbb{K}^{m\times n\times p}$
are 
\[
\left\{ \text{Prod}\left(\mathbf{P}^{\top^{2}},\mathbf{Q}^{\top^{2}},\text{Prod}\left(\text{Prod}\left(\mathbf{U},\mathbf{H},\mathbf{V}\right),\mathbf{E}^{\top},\mathbf{F}^{\top}\right)\right):\left(\mathbf{U},\mathbf{V}\right),\left(\mathbf{E},\mathbf{F}\right)\left(\mathbf{P},\mathbf{Q}\right)\in\text{GL}_{m\times n\times p}\left(m\times p\times p,\,p\times n\times p,\,\mathbb{K}\right)\right\} ,
\]
\[
\left\{ \text{Prod}\left(\mathbf{P}^{\top^{2}},\mathbf{Q}^{\top^{2}},\text{Prod}\left(\mathbf{U},\text{Prod}\left(\mathbf{H},\mathbf{E}^{\top},\mathbf{F}^{\top}\right),\mathbf{V}\right)\right):\left(\mathbf{U},\mathbf{V}\right),\left(\mathbf{E},\mathbf{F}\right)\left(\mathbf{P},\mathbf{Q}\right)\in\text{GL}_{m\times n\times p}\left(m\times p\times p,\,p\times n\times p,\,\mathbb{K}\right)\right\} ,
\]
\[
\left\{ \text{Prod}\left(\text{Prod}\left(\mathbf{P}^{\top^{2}},\mathbf{Q}^{\top^{2}},\text{Prod}\left(\mathbf{U},\mathbf{H},\mathbf{V}\right)\right),\mathbf{E}^{\top},\mathbf{F}^{\top}\right):\left(\mathbf{U},\mathbf{V}\right),\left(\mathbf{E},\mathbf{F}\right)\left(\mathbf{P},\mathbf{Q}\right)\in\text{GL}_{m\times n\times p}\left(m\times p\times p,\,p\times n\times p,\,\mathbb{K}\right)\right\} ,
\]
\[
\left\{ \text{Prod}\left(\text{Prod}\left(\mathbf{U},\text{Prod}\left(\mathbf{P}^{\top^{2}},\mathbf{Q}^{\top^{2}},\mathbf{H}\right),\mathbf{V}\right),\mathbf{E}^{\top},\mathbf{F}^{\top}\right):\left(\mathbf{U},\mathbf{V}\right),\left(\mathbf{E},\mathbf{F}\right)\left(\mathbf{P},\mathbf{Q}\right)\in\text{GL}_{m\times n\times p}\left(m\times p\times p,\,p\times n\times p,\,\mathbb{K}\right)\right\} ,
\]
\[
\left\{ \text{Prod}\left(\mathbf{U},\text{Prod}\left(\text{Prod}\left(\mathbf{P}^{\top^{2}},\mathbf{Q}^{\top^{2}},\mathbf{H}\right),\mathbf{E}^{\top},\mathbf{F}^{\top}\right),\mathbf{V}\right):\left(\mathbf{U},\mathbf{V}\right),\left(\mathbf{E},\mathbf{F}\right)\left(\mathbf{P},\mathbf{Q}\right)\in\text{GL}_{m\times n\times p}\left(m\times p\times p,\,p\times n\times p,\,\mathbb{K}\right)\right\} ,
\]
\[
\left\{ \text{Prod}\left(\mathbf{U},\text{Prod}\left(\mathbf{P}^{\top^{2}},\mathbf{Q}^{\top^{2}},\text{Prod}\left(\mathbf{H},\mathbf{E}^{\top},\mathbf{F}^{\top}\right)\right),\mathbf{V}\right):\left(\mathbf{U},\mathbf{V}\right),\left(\mathbf{E},\mathbf{F}\right)\left(\mathbf{P},\mathbf{Q}\right)\in\text{GL}_{m\times n\times p}\left(m\times p\times p,\,p\times n\times p,\,\mathbb{K}\right)\right\} .
\]
For convenience we adopt the notationanl convention such that 
\[
\text{Prod}\left(\mathbf{C},\text{Prod}\left(\mathbf{A},\mathbf{X},\mathbf{B}\right),\mathbf{D}\right)=\mathbf{X},\ \forall\,\mathbf{X}\in\mathbb{K}^{m\times n\times p}\Leftrightarrow\begin{cases}
\begin{array}{ccc}
\mathbf{C} & = & \mathbf{A}^{-1_{0}}\\
\\
\mathbf{D} & = & \mathbf{B}^{-1_{2}}
\end{array}\end{cases},
\]
\[
\text{Prod}\left(\text{Prod}\left(\mathbf{X}^{\top},\mathbf{B}^{\top},\mathbf{A}^{\top}\right),\mathbf{D}^{\top},\mathbf{C}^{\top}\right)=\mathbf{X},\ \forall\,\mathbf{X}\in\mathbb{K}^{m\times n\times p}\Leftrightarrow\begin{cases}
\begin{array}{ccc}
\mathbf{C}^{\top} & = & \left(\mathbf{A}^{\top}\right)^{-1_{2}}\\
\\
\mathbf{D}^{\top} & = & \left(\mathbf{B}^{\top}\right)^{-1_{1}}
\end{array}\end{cases},
\]
and
\[
\text{Prod}\left(\mathbf{D}^{\top^{2}},\mathbf{C}^{\top^{2}},\text{Prod}\left(\mathbf{B}^{\top^{2}},\mathbf{A}^{\top^{2}},\mathbf{X}^{\top^{2}}\right)\right)=\mathbf{X},\ \forall\,\mathbf{X}\in\mathbb{K}^{m\times n\times p}\Leftrightarrow\begin{cases}
\begin{array}{ccc}
\mathbf{C}^{\top^{2}} & = & \left(\mathbf{A}^{\top^{2}}\right)^{-1_{1}}\\
\\
\mathbf{D}^{\top^{2}} & = & \left(\mathbf{B}^{\top^{2}}\right)^{-1_{0}}
\end{array}\end{cases}.
\]

\section{SVD via Symmetrization.}

Recall the canonical $\mathbb{R}^{2\times2}$ representation of the
field $\mathbb{C}$ is prescribed by the correspondence 
\begin{equation}
\left(a+b\sqrt{-1}\right)\leftrightarrow\begin{pmatrix}a & -b\\
b & a
\end{pmatrix}.\label{representation_of_CC}
\end{equation}
We therefore express an arbitrary $\mathbf{M}\in\mathbb{C}^{n\times n}$
as a new matrix $\mathbf{M}^{\prime}\in\mathbb{R}^{2n\times2n}$ obtained
by replacing each entry of $\mathbf{M}$ by the corresponding $2\times2$
real matrix representation. It follows that no loss of generality
incurs from restricting the discussion to real matrices. 

It is well known that the Singular Value Decomposition (or SVD for
short ) of $\mathbf{A}\in\mathbb{R}^{n\times n}$ is obtained by solving
for matrices $\mathbf{U}$, $\mathbf{V}$, diag$\left(\boldsymbol{\mu}\right)$
and diag$\left(\boldsymbol{\nu}\right)$ in the constraints
\[
\begin{cases}
\begin{array}{ccc}
\left(\mathbf{A}\mathbf{A}^{\top}\right)^{k} & = & \left(\mathbf{U}\text{diag}\left(\boldsymbol{\mu}\right)^{k}\right)\left(\mathbf{U}\text{diag}\left(\boldsymbol{\mu}\right)^{k}\right)^{\top}\\
 & \text{and}\\
\left(\mathbf{A}^{\top}\mathbf{A}\right)^{k} & = & \left(\text{diag}\left(\boldsymbol{\nu}\right)^{k}\mathbf{V}\right)^{\top}\left(\text{diag}\left(\boldsymbol{\nu}\right)^{k}\mathbf{V}\right)
\end{array} & ,\forall\:k\in\left\{ 0,1\right\} \end{cases}.
\]
A distinctive feature of SVD constraints is that it can be equivalently
formulated as a pair of fixed point constraints of the form
\begin{equation}
\begin{cases}
\begin{array}{ccc}
\left(\mathbf{A}\mathbf{A}^{\top}\right)\left(\left(\mathbf{U}\text{diag}\left(\boldsymbol{\mu}\right)\right)^{\top}\right)^{-1} & = & \mathbf{U}\text{diag}\left(\boldsymbol{\mu}\right)\\
\\
\left(\left(\text{diag}\left(\boldsymbol{\nu}\right)\mathbf{V}\right)^{\top}\right)^{-1}\left(\mathbf{A}^{\top}\mathbf{A}\right) & = & \text{diag}\left(\boldsymbol{\nu}\right)\mathbf{V}
\end{array}.\end{cases}\label{Matrix Fixed Point Constraints}
\end{equation}
The fixed point formulation in Eq. (\ref{Matrix Fixed Point Constraints})
lies at the heart of iterative procedures for SVD numerical approximation
schemes which fortunately extend to hypermatrices. Characteristic
polynomials which eliminate the entries of $\mathbf{U}$ and $\mathbf{V}$
from the SVD constraints in Eq. (\ref{SVD constraints}) are 
\begin{equation}
\begin{cases}
\begin{array}{ccc}
\text{Rank}\left(\mathbf{A}\mathbf{A}^{\top}-\left(\mathbf{U}\mu_{i}\mathbf{I}_{n}\right)\left(\mathbf{U}\mu_{i}\mathbf{I}_{n}\right)^{\top}\right) & < & \text{Rank}\left(\mathbf{A}\mathbf{A}^{\top}\right)\\
 & \text{and}\\
\text{Rank}\left(\mathbf{A}^{\top}\mathbf{A}-\left(\nu_{i}\mathbf{I}_{n}\mathbf{V}\right)^{\top}\left(\nu_{i}\mathbf{I}_{n}\mathbf{V}\right)\right) & < & \text{Rank}\left(\mathbf{A}\mathbf{A}^{\top}\right)
\end{array}\end{cases}\implies\begin{cases}
\begin{array}{c}
\det\left(\mathbf{A}\mathbf{A}^{\top}-\mu_{i}^{2}\mathbf{I}_{n}\right)=0\\
\text{and}\\
\det\left(\mathbf{A}^{\top}\mathbf{A}-\nu_{i}^{2}\mathbf{I}_{n}\right)=0
\end{array}, & \forall\:0\le i<n\end{cases}
\end{equation}
It is well known that $\left\{ \mu_{i}^{2}:\,0\le i<n\right\} =\left\{ \nu_{i}^{2}:\,0\le i<n\right\} $,
and as a result we can take 
\[
\text{diag}\left(\boldsymbol{\mu}\right)=\text{diag}\left(\boldsymbol{\sigma}\right)=\text{diag}\left(\boldsymbol{\nu}\right).
\]
Once the singular values are known, we simultaneously solve for entries
of $\mathbf{U}$ via constraints given by
\[
\left(\mathbf{I}_{n}\otimes\text{Vandermonde}\left\{ \boldsymbol{\sigma}\circ\boldsymbol{\sigma}\right\} \right)\text{vec}\left(\U\left[i,k\right]\U\left[j,k\right]:\begin{array}{c}
0\le i<j<n\\
0\le k<n
\end{array}\right)=\text{vec}\left(\left(\mathbf{A}\mathbf{A}^{\top}\right)^{k}\left[i,j\right]:\begin{array}{c}
0\le i<j<n\\
0\le k<n
\end{array}\right)
\]
and also simultaneously solve for all entries of $\mathbf{V}$ via
constraints given by
\[
\left(\mathbf{I}_{n}\otimes\text{Vandermonde}\left\{ \boldsymbol{\sigma}\circ\boldsymbol{\sigma}\right\} \right)\text{vec}\left(\V\left[k,i\right]\V\left[k,j\right]:\begin{array}{c}
0\le i<j<n\\
0\le k<n
\end{array}\right)=\text{vec}\left(\left(\mathbf{A}^{\top}\mathbf{A}\right)^{k}\left[i,j\right]:\begin{array}{c}
0\le i<j<n\\
0\le k<n
\end{array}\right).
\]
Note that the constraints above express a composition of constraints
of type one and two as described in \cite{construct2018}.

We now extend to third order hypermatrices the matrix symmetrization
formulation of the SVD. For an arbitrary $\mathbf{A}\in\mathbb{C}^{n\times n\times n}$,
the three products of transposes which necessarily result in a symmetric
hypermatrix are
\[
\begin{cases}
\begin{array}{ccc}
\mbox{Prod}\left(\mathbf{A},\mathbf{A}^{\top^{2}},\mathbf{A}^{\top}\right)^{\top} & = & \mbox{Prod}\left(\mathbf{A},\mathbf{A}^{\top^{2}},\mathbf{A}^{\top}\right)\\
\\
\mbox{Prod}\left(\mathbf{A}^{\top},\mathbf{A},\mathbf{A}^{\top^{2}}\right)^{\top} & = & \mbox{Prod}\left(\mathbf{A}^{\top},\mathbf{A},\mathbf{A}^{\top^{2}}\right)\\
\\
\mbox{Prod}\left(\mathbf{A}^{\top^{2}},\mathbf{A}^{\top},\mathbf{\mathbf{A}}\right)^{\top} & = & \mbox{Prod}\left(\mathbf{A}^{\top^{2}},\mathbf{A}^{\top},\mathbf{\mathbf{A}}\right)
\end{array}.\end{cases}
\]
Just as was done for matrices, we devise the SVD from the spectral
decomposition of these symmetric products of transposes. Recall that
the scaling hypermatrices described in section 5 are hypermatrix analog
of diagonal matrices and characterized by the constraints
\[
\mathbf{D}^{\circ^{3}}\in\left\{ \text{Prod}\left(\mathbf{D}^{\top},\mathbf{D}^{\top^{2}},\mathbf{D}\right),\:\text{Prod}\left(\mathbf{D},\mathbf{D}^{\top},\mathbf{D}^{\top^{2}}\right),\:\text{Prod}\left(\mathbf{D}^{\top^{2}},\mathbf{D},\mathbf{D}^{\top}\right)\right\} ,
\]
where $\mathbf{D}^{\circ^{3}}$ represents the Hadamard exponent of
the scaling hypermatrix $\mathbf{D}$. Here we recall that the Hadamard
exponent $\mathbf{H}^{\circ^{z}}$ is defined for an arbitrary $\mathbf{H}\in\mathbb{C}^{m\times n\times p}$
and $z\in\mathbb{C}$ as follows

\[
\mathbf{H}^{\circ^{z}}\left[i,j,k\right]=\begin{cases}
\begin{array}{cc}
\left(\mathbf{H}\left[i,j,k\right]\right)^{z} & \text{ if }\mathbf{H}\left[i,j,k\right]\ne0\\
0 & \text{otherwise}
\end{array}.\end{cases}
\]
The above constraints are thus the hypermatrix diagonality constraints
generalized from the following matrix constraints 
\[
\text{Prod}\left(\mathbf{D}^{\top},\mathbf{D}\right)=\mathbf{D}^{\circ^{2}}=\text{Prod}\left(\mathbf{D},\mathbf{D}^{\top}\right).
\]
Note that in contrast to the matrix case, scaling hypermatrices are
not necessarily symmetric. For simplicity we describe the detailed
derivation of the SVD for an arbitrary side length two cubic hypermatrix
$\mathbf{A}\in\mathbb{C}^{2\times2\times2}$ whose entries are given
by
\[
\mathbf{A}\left[:,:,0\right]=\left(\begin{array}{rr}
a_{000} & a_{010}\\
a_{100} & a_{110}
\end{array}\right),\quad\mathbf{A}\left[:,:,1\right]=\left(\begin{array}{rr}
a_{001} & a_{011}\\
a_{101} & a_{111}
\end{array}\right),
\]
associated with the spectral decomposition constraints
\begin{equation}
\begin{cases}
\begin{array}{ccc}
\text{Prod}\left(\mathbf{A},\mathbf{A}^{\top^{2}},\mathbf{A}^{\top}\right) & = & \text{Prod}\left(\text{Prod}\left(\mathbf{U},\mathbf{D}_{\boldsymbol{\mu}},\mathbf{D}_{\boldsymbol{\mu}}^{\top}\right),\mbox{Prod}\left(\mathbf{U},\mathbf{D}_{\boldsymbol{\mu}},\mathbf{D}_{\boldsymbol{\mu}}^{\top}\right)^{\top^{2}},\mbox{Prod}\left(\mathbf{U},\mathbf{D}_{\boldsymbol{\mu}},\mathbf{D}_{\boldsymbol{\mu}}^{\top}\right)^{\top}\right)\\
\\
\mathbf{D}_{\boldsymbol{\mu}}^{\circ^{3}} & = & \text{Prod}\left(\mathbf{D}_{\boldsymbol{\mu}}^{\top},\mathbf{D}_{\boldsymbol{\mu}}^{\top^{2}},\mathbf{D}_{\boldsymbol{\mu}}\right)\\
\\
\mathbf{D}_{\boldsymbol{\mu}}\left[:,:,0\right] & = & \left(\begin{array}{rr}
\mu_{00} & 0\\
\mu_{01} & 0
\end{array}\right)\\
\\
\mathbf{D}_{\boldsymbol{\mu}}\left[:,:,1\right] & = & \left(\begin{array}{rr}
0 & \mu_{01}\\
0 & \mu_{11}
\end{array}\right)\\
\text{Prod}\left(\mathbf{U},\mathbf{D}_{\boldsymbol{\mu}},\mathbf{D}_{\boldsymbol{\mu}}^{\top}\right)\left[i,j,k\right] & = & \mu_{\min\left\{ i,j\right\} \max\left\{ i,j\right\} }\,\mu_{\min\left\{ j,k\right\} \max\left\{ j,k\right\} }\,u_{ijk}
\end{array},\end{cases}
\end{equation}
\begin{equation}
\begin{cases}
\begin{array}{ccc}
\text{Prod}\left(\mathbf{A}^{\top},\mathbf{A},\mathbf{A}^{\top^{2}}\right) & = & \text{Prod}\left(\text{Prod}\left(\mathbf{D}_{\boldsymbol{\nu}}^{\top},\mathbf{V},\mathbf{D}_{\boldsymbol{\nu}}\right)^{\top},\text{Prod}\left(\mathbf{D}_{\boldsymbol{\nu}}^{\top},\mathbf{V},\mathbf{D}_{\boldsymbol{\nu}}\right),\text{Prod}\left(\mathbf{D}_{\boldsymbol{\nu}}^{\top},\mathbf{V},\mathbf{D}_{\boldsymbol{\nu}}\right)^{\top^{2}}\right)\\
\\
\mathbf{D}_{\boldsymbol{\nu}}^{\circ^{3}} & = & \text{Prod}\left(\mathbf{D}_{\boldsymbol{\nu}},\mathbf{D}_{\boldsymbol{\nu}}^{\top},\mathbf{D}_{\boldsymbol{\nu}}^{\top^{2}}\right)\\
\\
\mathbf{D}_{\boldsymbol{\nu}}\left[:,:,0\right] & = & \left(\begin{array}{rr}
\nu_{00} & \nu_{01}\\
0 & 0
\end{array}\right)\\
\\
\mathbf{D}_{\boldsymbol{\nu}}\left[:,:,1\right] & = & \left(\begin{array}{rr}
0 & 0\\
\nu_{01} & \nu_{11}
\end{array}\right)\\
\text{Prod}\left(\mathbf{D}_{\boldsymbol{\nu}}^{\top},\mathbf{V},\mathbf{D}_{\boldsymbol{\nu}}\right)\left[i,j,k\right] & = & \nu_{\min\left\{ i,k\right\} \max\left\{ i,k\right\} }\,\nu_{\min\left\{ j,k\right\} \max\left\{ j,k\right\} }\,v_{ijk}
\end{array},\end{cases}
\end{equation}
and
\begin{equation}
\begin{cases}
\begin{array}{ccc}
\text{Prod}\left(\mathbf{A}^{\top^{2}},\mathbf{A}^{\top},\mathbf{A}\right) & = & \text{Prod}\left(\text{Prod}\left(\mathbf{D}_{\boldsymbol{\omega}},\mathbf{D}_{\boldsymbol{\omega}}^{\top},\mathbf{W}\right)^{\top^{2}},\text{Prod}\left(\mathbf{D}_{\boldsymbol{\omega}},\mathbf{D}_{\boldsymbol{\omega}}^{\top},\mathbf{W}\right)^{\top},\text{Prod}\left(\mathbf{D}_{\boldsymbol{\omega}},\mathbf{D}_{\boldsymbol{\omega}}^{\top},\mathbf{W}\right)\right)\\
\\
\mathbf{D}_{\boldsymbol{\omega}}^{\circ^{3}} & = & \text{Prod}\left(\mathbf{D}_{\boldsymbol{\omega}},\mathbf{D}_{\boldsymbol{\omega}}^{\top},\mathbf{D}_{\boldsymbol{\omega}}^{\top^{2}}\right)\\
\\
\mathbf{D}_{\boldsymbol{\omega}}\left[:,:,0\right] & = & \left(\begin{array}{rr}
\omega_{00} & 0\\
0 & \omega_{01}
\end{array}\right)\\
\\
\mathbf{D}_{\boldsymbol{\omega}}\left[:,:,1\right] & = & \left(\begin{array}{rr}
\omega_{01} & 0\\
0 & \omega_{11}
\end{array}\right)\\
\text{Prod}\left(\mathbf{D}_{\boldsymbol{\omega}},\mathbf{D}_{\boldsymbol{\omega}}^{\top},\mathbf{W}\right)\left[i,j,k\right] & = & \omega_{\min\left\{ i,k\right\} \max\left\{ i,k\right\} }\,\omega_{\min\left\{ i,j\right\} \max\left\{ i,j\right\} }\,w_{ijk}
\end{array}.\end{cases}
\end{equation}
The hypermatrices $\mathbf{U}$, $\mathbf{V}$ and $\mathbf{W}$ whose
individual slices correspond to eigenmatrices are subject to the following
third order orthogonality constraints
\begin{equation}
\mbox{Prod}\left(\mathbf{U},\mathbf{U}^{\top^{2}},\mathbf{U}^{\top}\right)\left[i,j,k\right]=\mbox{Prod}\left(\mathbf{V}^{\top},\mathbf{V},\mathbf{V}^{\top^{2}}\right)\left[i,j,k\right]=\mbox{Prod}\left(\mathbf{W}^{\top^{2}},\mathbf{W}^{\top},\mathbf{W}\right)\left[i,j,k\right]=\begin{cases}
\begin{array}{cc}
1 & \text{if }i=j=k\\
0 & \text{otherwise}
\end{array}.\end{cases}
\end{equation}
A distinctive feature of SVD constraints quite analogous to the matrix
setting is the equivalent formulation as fixed point constraints of
the form
\begin{equation}
\begin{cases}
\begin{array}{ccc}
\text{Prod}\left(\mathbf{U},\mathbf{D}_{\boldsymbol{\mu}},\mathbf{D}_{\boldsymbol{\mu}}^{\top}\right) & = & \text{Prod}\left(\text{Prod}\left(\mathbf{A}^{\top},\mathbf{A},\mathbf{A}^{\top^{2}}\right),\left(\mbox{Prod}\left(\mathbf{U},\mathbf{D}_{\boldsymbol{\mu}},\mathbf{D}_{\boldsymbol{\mu}}^{\top}\right)^{\top^{2}}\right)^{-1_{1}},\left(\mbox{Prod}\left(\mathbf{U},\mathbf{D}_{\boldsymbol{\mu}},\mathbf{D}_{\boldsymbol{\mu}}^{\top}\right)^{\top}\right)^{-1_{2}}\right),\\
\\
\text{Prod}\left(\mathbf{D}_{\boldsymbol{\nu}}^{\top},\mathbf{V},\mathbf{D}_{\boldsymbol{\nu}}\right) & = & \text{Prod}\left(\left(\text{Prod}\left(\mathbf{D}_{\boldsymbol{\nu}}^{\top},\mathbf{V},\mathbf{D}_{\boldsymbol{\nu}}\right)^{\top}\right)^{-1_{0}},\text{Prod}\left(\mathbf{A}^{\top},\mathbf{A},\mathbf{A}^{\top^{2}}\right),\left(\text{Prod}\left(\mathbf{D}_{\boldsymbol{\nu}}^{\top},\mathbf{V},\mathbf{D}_{\boldsymbol{\nu}}\right)^{\top^{2}}\right)^{-1_{2}}\right),\\
\\
\text{Prod}\left(\mathbf{D}_{\boldsymbol{\omega}},\mathbf{D}_{\boldsymbol{\omega}}^{\top},\mathbf{W}\right) & = & \text{Prod}\left(\left(\text{Prod}\left(\mathbf{D}_{\boldsymbol{\omega}},\mathbf{D}_{\boldsymbol{\omega}}^{\top},\mathbf{W}\right)^{\top^{2}}\right)^{-1_{0}},\left(\text{Prod}\left(\mathbf{D}_{\boldsymbol{\omega}},\mathbf{D}_{\boldsymbol{\omega}}^{\top},\mathbf{W}\right)^{\top}\right)^{-1_{1}},\text{Prod}\left(\mathbf{A}^{\top},\mathbf{A},\mathbf{A}^{\top^{2}}\right)\right).
\end{array}\end{cases}
\end{equation}
Just as in the matrix case, the characteristic polynomials which determine
the entries of the scaling hypermatrices (hypermatrix analog of the
singular values) are given by constraints of the form 
\[
\forall\,0\le i<2,\quad\text{Rank}\left(\mathbf{A},\mathbf{A}^{\top^{2}},\mathbf{A}^{\top}\right)>\text{Rank}\left\{ \text{Prod}\left(\mathbf{A},\mathbf{A}^{\top^{2}},\mathbf{A}^{\top}\right)-\text{Prod}\left(\widetilde{\mathbf{U}}_{i},\widetilde{\mathbf{U}}_{i}^{\top^{2}},\widetilde{\mathbf{U}}_{i}^{\top}\right)\right\} 
\]
\[
\text{where}
\]
\[
\widetilde{\mathbf{U}}_{i}=\text{Prod}\left(\mathbf{U},\mathbf{D}_{\boldsymbol{\mu}}^{\left[i\right]},\left(\mathbf{D}_{\boldsymbol{\mu}}^{\left[i\right]}\right)^{\top}\right).
\]
\[
\forall\,0\le i<2,\quad\text{Rank}\left(\mathbf{A}^{\top},\mathbf{A},\mathbf{A}^{\top^{2}}\right)>\text{Rank}\left\{ \text{Prod}\left(\mathbf{A}^{\top},\mathbf{A},\mathbf{A}^{\top^{2}}\right)-\text{Prod}\left(\widetilde{\mathbf{V}}_{i}^{\top},\widetilde{\mathbf{V}}_{i},\widetilde{\mathbf{V}}_{i}^{\top^{2}}\right)\right\} 
\]
\[
\text{where}
\]
\[
\widetilde{\mathbf{V}}_{i}=\text{Prod}\left(\left(\mathbf{D}_{\boldsymbol{\nu}}^{\left[i\right]}\right)^{\top},\mathbf{V},\mathbf{D}_{\boldsymbol{\nu}}^{\left[i\right]}\right).
\]
and
\[
\forall\,0\le i<2,\quad\text{Rank}\left(\mathbf{A}^{\top^{2}},\mathbf{A}^{\top},\mathbf{A}\right)>\text{Rank}\left\{ \text{Prod}\left(\mathbf{A}^{\top^{2}},\mathbf{A}^{\top},\mathbf{A}\right)-\text{Prod}\left(\widetilde{\mathbf{W}}_{i}^{\top^{2}},\widetilde{\mathbf{W}}_{i}^{\top},\widetilde{\mathbf{W}}_{i}\right)\right\} 
\]
\[
\text{where}
\]
\[
\widetilde{\mathbf{W}}_{i}=\text{Prod}\left(\mathbf{D}_{\boldsymbol{\omega}}^{\left[i\right]},\left(\mathbf{D}_{\boldsymbol{\omega}}^{\left[i\right]}\right)^{\top},\mathbf{W}\right).
\]
The entries of the scaling hypermatrices above are given by
\[
\begin{array}{ccc}
\mathbf{D}_{\boldsymbol{\mu}}^{\left[0\right]}\left[:,:,0\right] & = & \left(\begin{array}{rr}
\mu_{00} & 0\\
\mu_{01} & 0
\end{array}\right)\\
\\
\mathbf{D}_{\boldsymbol{\mu}}^{\left[0\right]}\left[:,:,1\right] & = & \left(\begin{array}{rr}
0 & \mu_{00}\\
0 & \mu_{01}
\end{array}\right)\\
\\
\\
\mathbf{D}_{\boldsymbol{\mu}}^{\left[1\right]}\left[:,:,0\right] & = & \left(\begin{array}{rr}
\mu_{01} & 0\\
\mu_{11} & 0
\end{array}\right)\\
\\
\mathbf{D}_{\boldsymbol{\mu}}^{\left[1\right]}\left[:,:,1\right] & = & \left(\begin{array}{rr}
0 & \mu_{01}\\
0 & \mu_{11}
\end{array}\right)
\end{array};\;\begin{array}{ccc}
\mathbf{D}_{\boldsymbol{\nu}}^{\left[0\right]}\left[:,:,0\right] & = & \left(\begin{array}{rr}
\nu_{00} & \nu_{01}\\
0 & 0
\end{array}\right)\\
\\
\mathbf{D}_{\boldsymbol{\nu}}^{\left[0\right]}\left[:,:,1\right] & = & \left(\begin{array}{rr}
0 & 0\\
\nu_{00} & \nu_{01}
\end{array}\right)\\
\\
\\
\mathbf{D}_{\boldsymbol{\nu}}^{\left[1\right]}\left[:,:,0\right] & = & \left(\begin{array}{rr}
\nu_{01} & \nu_{11}\\
0 & 0
\end{array}\right)\\
\\
\mathbf{D}_{\boldsymbol{\nu}}^{\left[1\right]}\left[:,:,1\right] & = & \left(\begin{array}{rr}
0 & 0\\
\nu_{01} & \nu_{11}
\end{array}\right)
\end{array};\;\begin{array}{ccc}
\mathbf{D}_{\boldsymbol{\omega}}^{\left[0\right]}\left[:,:,0\right] & = & \left(\begin{array}{rr}
\omega_{00} & 0\\
0 & \omega_{01}
\end{array}\right)\\
\\
\mathbf{D}_{\boldsymbol{\omega}}^{\left[0\right]}\left[:,:,1\right] & = & \left(\begin{array}{rr}
\omega_{00} & 0\\
0 & \omega_{01}
\end{array}\right)\\
\\
\\
\mathbf{D}_{\boldsymbol{\omega}}^{\left[1\right]}\left[:,:,0\right] & = & \left(\begin{array}{rr}
\omega_{01} & 0\\
0 & \omega_{11}
\end{array}\right)\\
\\
\mathbf{D}_{\boldsymbol{\omega}}^{\left[1\right]}\left[:,:,1\right] & = & \left(\begin{array}{rr}
\omega_{01} & 0\\
0 & \omega_{11}
\end{array}\right)
\end{array}.
\]
 Consequently characteristic polynomial constraints are expressed
by 
\begin{equation}
\forall\,0\le i<2,\quad\begin{cases}
\begin{array}{ccc}
0 & = & \det\left\{ \text{Prod}\left(\mathbf{A},\mathbf{A}^{\top^{2}},\mathbf{A}^{\top}\right)-\text{Prod}\left(\widetilde{\mathbf{U}}_{i},\widetilde{\mathbf{U}}_{i}^{\top^{2}},\widetilde{\mathbf{U}}_{i}^{\top}\right)\right\} ,\\
\\
0 & = & \det\left\{ \text{Prod}\left(\mathbf{A}^{\top},\mathbf{A},\mathbf{A}^{\top^{2}}\right)-\text{Prod}\left(\widetilde{\mathbf{V}}_{i}^{\top},\widetilde{\mathbf{V}}_{i},\widetilde{\mathbf{V}}_{i}^{\top^{2}}\right)\right\} ,\\
\\
0 & = & \det\left\{ \text{Prod}\left(\mathbf{A}^{\top^{2}},\mathbf{A}^{\top},\mathbf{A}\right)-\text{Prod}\left(\widetilde{\mathbf{W}}_{i}^{\top^{2}},\widetilde{\mathbf{W}}_{i}^{\top},\widetilde{\mathbf{W}}_{i}\right)\right\} .
\end{array}\end{cases}\label{Characteristic_polynomials}
\end{equation}
Using the hypermatrix determinant formula introduced by Gnang and
Yuval in \cite{Gnang2017238}, the corresponding constraints are expressed
as
\[
\begin{cases}
\begin{array}{c}
0=\left(\mu_{01}^{6}-a_{101}^{3}-a_{111}^{3}\right)\left(a_{000}a_{001}a_{100}+a_{010}a_{011}a_{110}\right)^{3}-\left(\mu_{00}^{6}-a_{000}^{3}-a_{010}^{3}\right)\left(a_{001}a_{100}a_{101}+a_{011}a_{110}a_{111}\right)^{3}\\
0=\left(\mu_{11}^{6}-a_{101}^{3}-a_{111}^{3}\right)\left(a_{000}a_{001}a_{100}+a_{010}a_{011}a_{110}\right)^{3}-\left(\mu_{01}^{6}-a_{000}^{3}-a_{010}^{3}\right)\left(a_{001}a_{100}a_{101}+a_{011}a_{110}a_{111}\right)^{3}\\
\\
0=\left(\nu_{01}^{6}-a_{110}^{3}-a_{111}^{3}\right)\left(a_{000}a_{010}a_{100}+a_{001}a_{011}a_{101}\right)^{3}-\left(\nu_{00}^{6}-a_{000}^{3}-a_{001}^{3}\right)\left(a_{010}a_{100}a_{110}+a_{011}a_{101}a_{111}\right)^{3}\\
0=\left(\nu_{11}^{6}-a_{110}^{3}-a_{111}^{3}\right)\left(a_{000}a_{010}a_{100}+a_{001}a_{011}a_{101}\right)^{3}-\left(\nu_{01}^{6}-a_{000}^{3}-a_{001}^{3}\right)\left(a_{010}a_{100}a_{110}+a_{011}a_{101}a_{111}\right)^{3}\\
\\
0=\left(\omega_{01}^{6}-a_{011}^{3}-a_{111}^{3}\right)\left(a_{000}a_{001}a_{010}+a_{100}a_{101}a_{110}\right)^{3}-\left(\omega_{00}^{6}-a_{000}^{3}-a_{100}^{3}\right)\left(a_{001}a_{010}a_{011}+a_{101}a_{110}a_{111}\right)^{3}\\
0=\left(\omega_{11}^{6}-a_{011}^{3}-a_{111}^{3}\right)\left(a_{000}a_{001}a_{010}+a_{100}a_{101}a_{110}\right)^{3}-\left(\omega_{01}^{6}-a_{000}^{3}-a_{100}^{3}\right)\left(a_{001}a_{010}a_{011}+a_{101}a_{110}a_{111}\right)^{3}
\end{array}\end{cases}
\]
Once we have determined the entries of the scaling values, we simultaneously
solve for all entries of $\mathbf{U}$ via constraints given by
\[
\left(\begin{array}{rrrrrrrr}
1 & 1 & 0 & 0 & 0 & 0 & 0 & 0\\
\mu_{00}^{6} & \mu_{01}^{6} & 0 & 0 & 0 & 0 & 0 & 0\\
0 & 0 & 1 & 1 & 0 & 0 & 0 & 0\\
0 & 0 & \mu_{00}^{4}\mu_{01}^{2} & \mu_{01}^{4}\mu_{11}^{2} & 0 & 0 & 0 & 0\\
0 & 0 & 0 & 0 & 1 & 1 & 0 & 0\\
0 & 0 & 0 & 0 & \mu_{00}^{2}\mu_{01}^{4} & \mu_{01}^{2}\mu_{11}^{4} & 0 & 0\\
0 & 0 & 0 & 0 & 0 & 0 & 1 & 1\\
0 & 0 & 0 & 0 & 0 & 0 & \mu_{01}^{6} & \mu_{11}^{6}
\end{array}\right)\cdot\left(\begin{array}{c}
u_{000}^{3}\\
u_{010}^{3}\\
u_{000}u_{001}u_{100}\\
u_{010}u_{011}u_{110}\\
u_{001}u_{100}u_{101}\\
u_{011}u_{110}u_{111}\\
u_{101}^{3}\\
u_{111}^{3}
\end{array}\right)=\left(\begin{array}{c}
1\\
a_{000}^{3}+a_{010}^{3}\\
0\\
a_{000}a_{001}a_{100}+a_{010}a_{011}a_{110}\\
0\\
a_{001}a_{100}a_{101}+a_{011}a_{110}a_{111}\\
1\\
a_{101}^{3}+a_{111}^{3}
\end{array}\right),
\]
solve for all entries of $\mathbf{V}$ via constraints given by
\[
\left(\begin{array}{rrrrrrrr}
1 & 1 & 0 & 0 & 0 & 0 & 0 & 0\\
\nu_{00}^{6} & \nu_{01}^{6} & 0 & 0 & 0 & 0 & 0 & 0\\
0 & 0 & 1 & 1 & 0 & 0 & 0 & 0\\
0 & 0 & \nu_{00}^{4}\nu_{01}^{2} & \nu_{01}^{4}\nu_{11}^{2} & 0 & 0 & 0 & 0\\
0 & 0 & 0 & 0 & 1 & 1 & 0 & 0\\
0 & 0 & 0 & 0 & \nu_{00}^{2}\nu_{01}^{4} & \nu_{01}^{2}\nu_{11}^{4} & 0 & 0\\
0 & 0 & 0 & 0 & 0 & 0 & 1 & 1\\
0 & 0 & 0 & 0 & 0 & 0 & \nu_{01}^{6} & \nu_{11}^{6}
\end{array}\right)\cdot\left(\begin{array}{c}
v_{000}^{3}\\
v_{001}^{3}\\
v_{000}v_{010}v_{100}\\
v_{001}v_{011}v_{101}\\
v_{010}v_{100}v_{110}\\
v_{011}v_{101}v_{111}\\
v_{110}^{3}\\
v_{111}^{3}
\end{array}\right)=\left(\begin{array}{c}
1\\
a_{000}^{3}+a_{001}^{3}\\
0\\
a_{000}a_{010}a_{100}+a_{001}a_{011}a_{101}\\
0\\
a_{010}a_{100}a_{110}+a_{011}a_{101}a_{111}\\
1\\
a_{110}^{3}+a_{111}^{3}
\end{array}\right),
\]
and also solve for all entries of $\mathbf{W}$ via constraints given
by
\[
\left(\begin{array}{rrrrrrrr}
1 & 1 & 0 & 0 & 0 & 0 & 0 & 0\\
\omega_{00}^{6} & \omega_{01}^{6} & 0 & 0 & 0 & 0 & 0 & 0\\
0 & 0 & 1 & 1 & 0 & 0 & 0 & 0\\
0 & 0 & \omega_{00}^{4}\omega_{01}^{2} & \omega_{01}^{4}\omega_{11}^{2} & 0 & 0 & 0 & 0\\
0 & 0 & 0 & 0 & 1 & 1 & 0 & 0\\
0 & 0 & 0 & 0 & \omega_{00}^{2}\omega_{01}^{4} & \omega_{01}^{2}\omega_{11}^{4} & 0 & 0\\
0 & 0 & 0 & 0 & 0 & 0 & 1 & 1\\
0 & 0 & 0 & 0 & 0 & 0 & \omega_{01}^{6} & \omega_{11}^{6}
\end{array}\right)\cdot\left(\begin{array}{c}
w_{000}^{3}\\
w_{100}^{3}\\
w_{000}w_{001}w_{010}\\
w_{100}w_{101}w_{110}\\
w_{001}w_{010}w_{011}\\
w_{101}w_{110}w_{111}\\
w_{011}^{3}\\
w_{111}^{3}
\end{array}\right)=\left(\begin{array}{c}
1\\
a_{000}^{3}+a_{100}^{3}\\
0\\
a_{000}a_{001}a_{010}+a_{100}a_{101}a_{110}\\
0\\
a_{001}a_{010}a_{011}+a_{101}a_{110}a_{111}\\
1\\
a_{011}^{3}+a_{111}^{3}
\end{array}\right)
\]
Note that constraints above correspond to a composition of constraints
of type one and two discussed in \cite{construct2018}. 

The hypermatrix SVD is thus expressed by the following sum of outer
products 
\begin{equation}
\mathbf{A}=\sum_{0\le i,j,k<2}\sigma_{i,j,k}\,\text{Prod}\left(\widetilde{\mathbf{U}}\left[:,i,:\right],\,\widetilde{\mathbf{V}}\left[:,:,j\right],\,\widetilde{\mathbf{W}}\left[k,:,:\right]\right)\label{SVD}
\end{equation}
\[
\text{where}
\]
\[
\widetilde{\mathbf{U}}=\text{Prod}\left(\mathbf{U},\mathbf{D}_{\boldsymbol{\mu}},\mathbf{D}_{\boldsymbol{\mu}}^{\top}\right),\quad\,\widetilde{\mathbf{V}}=\text{Prod}\left(\mathbf{D}_{\boldsymbol{\nu}}^{\top},\mathbf{V},\mathbf{D}_{\boldsymbol{\nu}}\right),\:\text{ and }\:\widetilde{\mathbf{W}}=\text{Prod}\left(\mathbf{D}_{\boldsymbol{\omega}},\mathbf{D}_{\boldsymbol{\omega}}^{\top},\mathbf{W}\right).
\]
The coefficients $\left\{ \sigma_{i,j,k}\,:\,0\le i,j,k<2\right\} \subset\mathbb{C}$
of the linear combination in Eq. (\ref{SVD}) are obtained through
solving a system of linear equations. The expansion in Eq. (\ref{SVD})
is equivalently expressed as
\[
\mathbf{A}=\text{Prod}\left(\mathbf{U}^{\prime},\,\mathbf{V}^{\prime},\,\mathbf{W}^{\prime}\right),
\]
where $\mathbf{U}^{\prime}\in\mathbb{C}^{2\times\left\Vert \boldsymbol{\sigma}\right\Vert _{\ell_{0}}\times2}$,
$\mathbf{V}^{\prime}\in\mathbb{C}^{2\times2\times\left\Vert \boldsymbol{\sigma}\right\Vert _{\ell_{0}}}$
, and $\mathbf{W}^{\prime}\in\mathbb{C}^{\left\Vert \boldsymbol{\sigma}\right\Vert _{\ell_{0}}\times2\times2}$,
and $\boldsymbol{\sigma}$ is the vector whose entries are made up
of the coefficients $\left\{ \sigma_{i,j,k}\,:\,0\le i,j,k<2\right\} $
in the linear combination. As an illustration, consider the task of
expressing the SVD of hypermatrices of arbitrary side lengths generated
from $2\times2\times2$ hypermatrices by taking combinations of direct
sums and Kronecker products. As shown in \cite{gnang2018spectral},
and similarly to the matrix case, when given the SVD of hypermatrices
$\mathbf{A}_{0}\in\mathbb{C}^{m\times m\times m}$ and $\mathbf{A}_{1}\in\mathbb{C}^{n\times n\times n}$
\[
\mathbf{A}_{0}=\text{Prod}\left(\mathbf{U}_{0}^{\prime},\,\mathbf{V}_{0}^{\prime},\,\mathbf{W}_{0}^{\prime}\right),\quad\mathbf{A}_{1}=\text{Prod}\left(\mathbf{U}_{1}^{\prime},\,\mathbf{V}_{1}^{\prime},\,\mathbf{W}_{1}^{\prime}\right),
\]
then the SVD of $\mathbf{A}_{0}\otimes\mathbf{A}_{1}$ and $\mathbf{A}_{0}\oplus\mathbf{A}_{1}$
are expressed by
\[
\mathbf{A}_{0}\otimes\mathbf{A}_{1}=\text{Prod}\left(\mathbf{U}_{0}^{\prime}\otimes\mathbf{U}_{1}^{\prime},\,\mathbf{V}_{0}^{\prime}\otimes\mathbf{V}_{1}^{\prime},\,\mathbf{W}_{0}^{\prime}\otimes\mathbf{W}_{1}^{\prime}\right)\in\mathbb{C}^{mn\times mn\times mn},
\]
\[
\mathbf{A}_{0}\oplus\mathbf{A}_{1}=\text{Prod}\left(\mathbf{U}_{0}^{\prime}\oplus\mathbf{U}_{1}^{\prime},\,\mathbf{V}_{0}^{\prime}\oplus\mathbf{V}_{1}^{\prime},\,\mathbf{W}_{0}^{\prime}\oplus\mathbf{W}_{1}^{\prime}\right)\in\mathbb{C}^{\left(m+n\right)\times\left(m+n\right)\times\left(m+n\right)},
\]

\section{Action on vector spaces and orthogonality.}

\subsection{The matrix case.}

The action of a matrix in $\mathbb{C}^{n\times n}$ on the vector
space $\mathbb{C}^{n\times1}$ can be seen as a special instance of
a more general (not necessarily linear) map introduced in \cite{Gnang2017238}
specified in terms of a matrix pair $\left(\mathbf{A},\mathbf{B}\right)\in\mathbb{C}^{n\times n}\times\mathbb{C}^{n\times n}$
as follows
\[
\mathcal{T}_{\mathbf{A},\mathbf{B}}\,:\,\mathbb{C}^{n\times1}\rightarrow\mathbb{C}^{n\times1},\quad\mathbf{y}=\mathcal{T}_{\mathbf{A},\mathbf{B}}\left(\mathbf{x}\right),
\]
\[
\text{such that}
\]
\begin{equation}
\forall\:0\le k<n,\quad\mathbf{y}\left[k\right]=\sqrt{\text{Prod}_{\mathbf{P}_{k}}\left(\mathbf{x}^{\top},\mathbf{x}\right)}\ \text{ where }\ \mathbf{P}_{k}=\text{Prod}_{\mathbf{I}_{n}\left[:,k\right]\mathbf{I}_{n}\left[k,:\right]}\left(\mathbf{A},\mathbf{B}\right).
\end{equation}
Note that the map $\mathcal{T}_{\mathbf{A},\mathbf{B}}$ is determined
up to the sign of the entries of its output. Invertibility in this
context means that neither of the $n$ univariate polynomials in 
\[
\text{Resultant}_{\mathbf{x}}\left\{ \text{Prod}_{\mathbf{P}_{k}}\left(\mathbf{x}^{\top},\mathbf{x}\right):0\le k<n\right\} ,
\]
is an identically non-zero constant. For instance when $n=2$ and
\[
\mathbf{A}=\left(\begin{array}{rr}
a_{00} & a_{01}\\
a_{10} & a_{11}
\end{array}\right),\;\mathbf{B}=\left(\begin{array}{rr}
b_{00} & b_{01}\\
b_{10} & b_{11}
\end{array}\right),
\]
the map $\mathcal{T}_{\mathbf{A},\mathbf{B}}$ is invertible if neither
of the polynomials in
\[
\left\{ Q_{0}\left(x_{0}\right),\,Q_{1}\left(x_{1}\right)\right\} =\text{Resultant}_{\mathbf{x}}\left\{ \text{Prod}_{\mathbf{P}_{k}}\left(\mathbf{x}^{\top},\mathbf{x}\right):0\le k<2\right\} 
\]
explicitly given by 
\[
Q_{0}\left(x_{0}\right)=\left(a_{01}b_{10}x_{0}^{2}-y_{1}^{2}\right)^{2}a_{10}^{2}b_{01}^{2}-2\,\left(a_{00}b_{00}x_{0}^{2}-y_{0}^{2}\right)\left(a_{01}b_{10}x_{0}^{2}-y_{1}^{2}\right)a_{10}a_{11}b_{01}b_{11}+\left(a_{00}b_{00}x_{0}^{2}-y_{0}^{2}\right)^{2}a_{11}^{2}b_{11}^{2}+
\]
\[
\left(-1\right)\left(a_{01}b_{10}x_{0}^{2}-y_{1}^{2}\right)\left(a_{10}b_{00}x_{0}+a_{00}b_{01}x_{0}\right)\left(a_{11}b_{10}x_{0}+a_{01}b_{11}x_{0}\right)a_{10}b_{01}+\left(a_{00}b_{00}x_{0}^{2}-y_{0}^{2}\right)\left(a_{11}b_{10}x_{0}+a_{01}b_{11}x_{0}\right)^{2}a_{10}b_{01}+
\]
\[
\left(a_{01}b_{10}x_{0}^{2}-y_{1}^{2}\right)\left(a_{10}b_{00}x_{0}+a_{00}b_{01}x_{0}\right)^{2}a_{11}b_{11}-\left(a_{00}b_{00}x_{0}^{2}-y_{0}^{2}\right)\left(a_{10}b_{00}x_{0}+a_{00}b_{01}x_{0}\right)\left(a_{11}b_{10}x_{0}+a_{01}b_{11}x_{0}\right)a_{11}b_{11}.
\]
and
\[
Q_{1}\left(x_{1}\right)=\left(a_{11}b_{11}x_{1}^{2}-y_{1}^{2}\right)^{2}a_{00}^{2}b_{00}^{2}-2\,\left(a_{10}b_{01}x_{1}^{2}-y_{0}^{2}\right)\left(a_{11}b_{11}x_{1}^{2}-y_{1}^{2}\right)a_{00}a_{01}b_{00}b_{10}+\left(a_{10}b_{01}x_{1}^{2}-y_{0}^{2}\right)^{2}a_{01}^{2}b_{10}^{2}+
\]
\[
\left(-1\right)\left(a_{11}b_{11}x_{1}^{2}-y_{1}^{2}\right)\left(a_{10}b_{00}x_{1}+a_{00}b_{01}x_{1}\right)\left(a_{11}b_{10}x_{1}+a_{01}b_{11}x_{1}\right)a_{00}b_{00}+\left(a_{10}b_{01}x_{1}^{2}-y_{0}^{2}\right)\left(a_{11}b_{10}x_{1}+a_{01}b_{11}x_{1}\right)^{2}a_{00}b_{00}+
\]
\[
\left(a_{11}b_{11}x_{1}^{2}-y_{1}^{2}\right)\left(a_{10}b_{00}x_{1}+a_{00}b_{01}x_{1}\right)^{2}a_{01}b_{10}-\left(a_{10}b_{01}x_{1}^{2}-y_{0}^{2}\right)\left(a_{10}b_{00}x_{1}+a_{00}b_{01}x_{1}\right)\left(a_{11}b_{10}x_{1}+a_{01}b_{11}x_{1}\right)a_{01}b_{10}.
\]
is an identically non-zero constant. Furthermore when $\mathbf{A}=\mathbf{B}^{-1}$
the map $\mathcal{T}_{\mathbf{A},\mathbf{B}}$ is subject to the resolution
of identity 
\[
\mathbf{y}=\mathcal{T}_{\mathbf{A},\mathbf{B}}\left(\mathbf{x}\right)\implies\text{Prod}\left(\mathbf{y}^{\top},\mathbf{y}\right)=\text{Prod}\left(\mathbf{x}^{\top},\mathbf{x}\right).
\]
In other words the map preserves the sum of squares of the entries.
Also note that when $\mathbf{A}=\mathbf{B}^{\top}$, the map $\mathcal{T}_{\mathbf{A},\mathbf{B}}$
expresses up to the sign of the entries a linear transformation. In
particular, when $\mathbf{A}\mathbf{B}=\mathbf{I}_{n}$ and $\mathbf{B}=\mathbf{A}^{\top}\in\mathbb{R}^{n\times n}$
the map $\mathcal{T}_{\mathbf{A},\mathbf{A}^{\top}}$ expresses up
to the entry signs a linear isometry of $\mathbb{R}^{n\times1}$,
thereby emphasizing the importance of matrix orthogonality. Recall
for illustration purposes that
\[
\mathbf{X}=\left(\begin{array}{cc}
x_{0} & x_{2}\\
x_{1} & x_{3}
\end{array}\right),
\]
is orthogonal if $\mathbf{X}\cdot\mathbf{X}^{\top}=\mathbf{I}_{2}$.
Hence 
\[
\mathbf{X}\cdot\mathbf{X}^{\top}=\mathbf{I}_{2}\implies\begin{cases}
\begin{array}{ccccc}
x_{0}^{2} & + & x_{2}^{2} & = & 1\\
x_{0}x_{1} & + & x_{2}x_{3} & = & 0\\
x_{1}^{2} & + & x_{3}^{2} & = & 1
\end{array}.\end{cases}
\]
On the one hand,
\[
0=\left(\prod_{0\le i<4}x_{i}\right)\implies\mathbf{X}\in\left\{ \left(\begin{array}{cc}
1 & 0\\
0 & 1
\end{array}\right),\,\left(\begin{array}{cc}
-1 & 0\\
0 & 1
\end{array}\right),\,\left(\begin{array}{cc}
1 & 0\\
0 & -1
\end{array}\right),\,\left(\begin{array}{cc}
-1 & 0\\
0 & -1
\end{array}\right)\right\} .
\]
On the other hand, when $0\ne\left(\underset{0\le i<4}{\prod}x_{i}\right)$
implies that 
\[
0=x_{0}x_{1}+x_{2}x_{3}\Leftrightarrow x_{0}x_{1}x_{2}^{-1}x_{3}^{-1}=-1,\ \forall\,k\in\mathbb{Z}.
\]
\begin{equation}
\implies\left(\begin{array}{c}
x_{0}\\
x_{1}\\
x_{2}\\
x_{3}
\end{array}\right)=\left(\begin{array}{c}
\nicefrac{-s\,t}{r}\\
r\\
s\\
t
\end{array}\right)\implies\mathbf{X}=\left(\begin{array}{cc}
\nicefrac{-s\,t}{r} & s\\
r & t
\end{array}\right)
\end{equation}
By normalizing the row of $\mathbf{X}$ we obtain the following parametrization
of the orthogonal matrices 
\begin{equation}
\mathbf{X}=\left(\begin{array}{ccc}
\frac{\nicefrac{-s\,t}{r}}{\sqrt{\left(\nicefrac{-st}{r}\right)^{2}+s^{2}}} &  & \frac{s}{\sqrt{\left(\nicefrac{-st}{r}\right)^{2}+s^{2}}}\\
\\
\frac{r}{\sqrt{r^{2}+t^{2}}} &  & \frac{t}{\sqrt{r^{2}+t^{2}}}
\end{array}\right),\;s\in\left\{ -1,1\right\} ,\ \text{ and }\ r\ne t\sqrt{-1}.\label{Orthogonal Matrix Parametrization}
\end{equation}
To express some important invariants of orthogonal matrices, consider
the index rotation operation introduced in \cite{gnang2018spectral},
noted $\mathbf{A}^{\text{R}_{\theta}}$ for $\theta\in\left\{ 0,\frac{\pi}{2},\pi,\frac{3\pi}{2}\right\} $,
which generalizes the matrix transpose operation and is defined for
an arbitrary $\mathbf{A}\in\mathbb{C}^{n\times n}$ as
\[
\A^{R_{0}}=\A,\ \A^{R_{\frac{\pi}{2}}}=\A^{\top}\mathbf{Q},\ \A^{R_{\pi}}=\mathbf{Q}\A\mathbf{Q},\ \A^{R_{3\pi/2}}=\mathbf{Q}\A^{\top},
\]
\[
\text{where}
\]
\[
\mathbf{Q}=\sum_{0\le i<n}\mathbf{I}_{n}\left[:,n-i-1\right]\mathbf{I}_{n}\left[i,:\right].
\]
Alternatively, we can also express the index rotation operation entry-wise
as 
\[
\left(\mathbf{A}^{\text{R}_{\theta}}\right)\left[i,j\right]=\mathbf{A}\left[\left(i-\frac{n-1}{2}\right)\cos\theta+\left(\frac{n-1}{2}-j\right)\sin\theta+\frac{n-1}{2},\left(i-\frac{n-1}{2}\right)\sin\theta-\left(\frac{n-1}{2}-j\right)\cos\theta+\frac{n-1}{2}\right].
\]
For instance for a given $3\times3$ matrix $\A$, we have
\[
\mathbf{A}^{\text{R}_{0}}=\left(\begin{array}{rrr}
a_{00} & a_{01} & a_{02}\\
a_{10} & a_{11} & a_{12}\\
a_{20} & a_{21} & a_{22}
\end{array}\right),\,\mathbf{A}^{\text{R}_{\frac{\pi}{2}}}=\left(\begin{array}{rrr}
a_{20} & a_{10} & a_{00}\\
a_{21} & a_{11} & a_{01}\\
a_{22} & a_{12} & a_{02}
\end{array}\right),\,\mathbf{A}^{\text{R}_{\pi}}=\left(\begin{array}{rrr}
a_{22} & a_{21} & a_{20}\\
a_{12} & a_{11} & a_{10}\\
a_{02} & a_{01} & a_{00}
\end{array}\right),\,\mathbf{A}^{\text{R}_{\frac{3\pi}{2}}}=\left(\begin{array}{rrr}
a_{02} & a_{12} & a_{22}\\
a_{01} & a_{11} & a_{21}\\
a_{00} & a_{10} & a_{20}
\end{array}\right).
\]
Following immediately from the orthogonal matrix parametrization in
Eq. (\ref{Orthogonal Matrix Parametrization}), we can obtain the
properties of orthogonal matrices,
\[
\mathbf{X}\mathbf{X}^{\top}=\mathbf{I}_{2}\implies\mathbf{X}^{\text{R}_{\frac{2\pi k}{4}}}\left(\mathbf{X}^{\text{R}_{\frac{2\pi k}{4}}}\right)^{\top}=\mathbf{I}_{2},\quad\forall\,k\in\left[0,4\right)\cap\mathbb{Z}
\]
\[
\text{and}
\]
\[
\mathbf{X}\mathbf{X}^{\top}=\mathbf{I}_{2}\implies\left(\mathbf{X}^{\text{R}_{\frac{2\pi k}{4}}}\right)^{\top}\mathbf{X}^{\text{R}_{\frac{2\pi k}{4}}}=\mathbf{I}_{2},\quad\forall\,k\in\left[0,4\right)\cap\mathbb{Z}.
\]
Furthermore, given $\mathbf{X}\mathbf{X}^{\top}=\mathbf{I}_{2}=$
$\mathbf{Y}\mathbf{Y}^{\top}$ we have 
\[
\left(\mathbf{X}\mathbf{Y}\right)\left(\mathbf{X}\mathbf{Y}\right)^{\top}=\mathbf{I}_{2},\:\left(\mathbf{X}\oplus\mathbf{Y}\right)\left(\mathbf{X}\oplus\mathbf{Y}\right)^{\top}=\mathbf{I}_{2}\oplus\mathbf{I}_{2},\:\text{and}\ \left(\mathbf{X}\otimes\mathbf{Y}\right)\left(\mathbf{X}\otimes\mathbf{Y}\right)^{\top}=\mathbf{I}_{2}\otimes\mathbf{I}_{2}.
\]
The canonical matrix representation of complex number described in
Eq. (\ref{representation_of_CC}) motivates a variant of the transpose
and index rotation operation which operates block partitioned matrices.
More precisely, consider the variant of the transpose and index rotation
operations defined on block matrices where each block is a square
matrix of the same size 

\[
\left(\begin{array}{cc}
\mathbf{A}_{00} & \mathbf{A}_{01}\\
\mathbf{A}_{10} & \mathbf{A}_{11}
\end{array}\right)^{\top_{b}}=\sum_{0\le i,j<2}\left(\mathbf{I}_{2}\left[:,i\right]\mathbf{I}_{2}\left[j,:\right]\right)^{\top}\otimes\mathbf{A}_{ij},\quad\left(\begin{array}{cc}
\mathbf{A}_{00} & \mathbf{A}_{01}\\
\mathbf{A}_{10} & \mathbf{A}_{11}
\end{array}\right)^{\top_{e}}=\sum_{0\le i,j<2}\left(\mathbf{I}_{2}\left[:,i\right]\mathbf{I}_{2}\left[j,:\right]\right)\otimes\mathbf{A}_{ij}^{\top},
\]
where $\left\{ \mathbf{A}_{00},\,\mathbf{A}_{01},\,\mathbf{A}_{10},\,\mathbf{A}_{11}\right\} $
correspond square matrix blocks all of the same size. Similarly 
\[
\left(\begin{array}{cc}
\mathbf{A}_{00} & \mathbf{A}_{01}\\
\mathbf{A}_{10} & \mathbf{A}_{11}
\end{array}\right)^{\text{R}_{_{b},\frac{2\pi k}{4}}}=\sum_{0\le i,j<2}\left(\mathbf{I}_{2}\left[:,i\right]\cdot\mathbf{I}_{2}\left[j,:\right]\right)^{\text{R}_{\frac{2\pi k}{4}}}\otimes\mathbf{A}_{ij},\quad\left(\begin{array}{cc}
\mathbf{A}_{00} & \mathbf{A}_{01}\\
\mathbf{A}_{10} & \mathbf{A}_{11}
\end{array}\right)^{\text{R}_{e,\frac{2\pi k}{4}}}=\sum_{0\le i,j<2}\left(\mathbf{I}_{2}\left[:,i\right]\cdot\mathbf{I}_{2}\left[j,:\right]\right)\otimes\mathbf{A}_{ij}^{\text{R}_{\frac{2\pi k}{4}}}.
\]
These operations distinguish actions deinfed on individual block matrices
from action defined on the whole matrix. This distinction will enable
us to generalize the matrix conjugate transpose operation. Note that
if $\mathbf{A}\in\left(O_{m}\left(\mathbb{C}\right)\right)^{n\times n}$
then 
\[
\forall\:0\le i<n,\quad\left(\frac{\mathbf{A}}{\sqrt{n}}\left(\frac{\mathbf{A}}{\sqrt{n}}^{\top_{e}}\right)^{\top_{b}}\right)\left[i,i\right]=\mathbf{I}_{m}=\left(\left(\frac{\mathbf{A}}{\sqrt{n}}^{\top_{e}}\right)^{\top_{b}}\frac{\mathbf{A}}{\sqrt{n}}\right)\left[i,i\right]
\]
It therefore follows that 
\[
\frac{\mathbf{A}}{\sqrt{n}}\left(\frac{\mathbf{A}}{\sqrt{n}}^{\top_{e}}\right)^{\top_{b}}=\mathbf{0}_{n\times n}\otimes\mathbf{0}_{m\times m}\implies\mathbf{A}=\mathbf{0}_{n\times n}\otimes\mathbf{0}_{m\times m}.
\]
The non-negativity property still holds if each block entry of $\mathbf{A}$
is positive scaling on an orthogonal matrix i. e. 
\[
\forall\:0\le i,j<n,\quad\mathbf{A}\left[i,j\right]=r_{ij}\,\mathbf{A}_{ij},\ \text{ where }r_{ij}>0\;\text{ and }\mathbf{A}_{ij}\in O_{m}\left(\mathbb{C}\right).
\]
In particular $\mathbf{A}\in\left(\mathbb{C}^{m\times m}\right)^{n\times n}$
such that 
\[
\forall\:0\le i,j<n,\quad\mathbf{A}\left[i,j\right]=r_{ij}\,\mathbf{A}_{ij},\ \text{ where }r_{ij}>0\;\text{ and }\mathbf{A}_{ij}\in O_{m}\left(\mathbb{C}\right)
\]
is called \emph{block unitary} if 
\[
\frac{\mathbf{A}}{\sqrt{n}}\left(\frac{\mathbf{A}}{\sqrt{n}}^{\top_{e}}\right)^{\top_{b}}=\mathbf{I}_{n}\otimes\mathbf{I}_{m}=\left(\frac{\mathbf{A}}{\sqrt{n}}^{\top_{e}}\right)^{\top_{b}}\frac{\mathbf{A}}{\sqrt{n}}
\]
In the case of $2\times2$ block matrices, block unitary constrains
for matrix blocks 
\[
\left\{ \mathbf{A}_{00},\,\mathbf{A}_{01},\,\mathbf{A}_{10},\,\mathbf{A}_{11}\right\} \subset O_{m}\left(\mathbb{C}\right)
\]
are expressed by 
\[
\left(\begin{array}{ccc}
\frac{\mathbf{A}_{00}}{\sqrt{2}} &  & \frac{\mathbf{A}_{01}}{\sqrt{2}}\\
\\
\frac{\mathbf{A}_{10}}{\sqrt{2}} &  & \frac{\mathbf{A}_{11}}{\sqrt{2}}
\end{array}\right)\left(\begin{array}{ccc}
\frac{\mathbf{A}_{00}^{\top}}{\sqrt{2}} &  & \frac{\mathbf{A}_{10}^{\top}}{\sqrt{2}}\\
\\
\frac{\mathbf{A}_{01}^{\top}}{\sqrt{2}} &  & \frac{\mathbf{A}_{11}^{\top}}{\sqrt{2}}
\end{array}\right)=\left(\begin{array}{ccc}
\mathbf{I}_{m} &  & \frac{\mathbf{A}_{00}\mathbf{A}_{10}^{\top}+\mathbf{A}_{01}\mathbf{A}_{11}^{\top}}{2}\\
\\
\frac{\mathbf{A}_{10}\mathbf{A}_{00}^{\top}+\mathbf{A}_{11}\mathbf{A}_{01}^{\top}}{2} &  & \mathbf{I}_{m}
\end{array}\right)
\]
which yields the constraints 
\[
\begin{cases}
\begin{array}{ccc}
\mathbf{A}_{00}\mathbf{A}_{10}^{\top}+\mathbf{A}_{01}\mathbf{A}_{11}^{\top} & = & \mathbf{0}_{m\times m}\\
\\
\mathbf{A}_{10}\mathbf{A}_{00}^{\top}+\mathbf{A}_{11}\mathbf{A}_{01}^{\top} & = & \mathbf{0}_{m\times m}
\end{array} & \implies\end{cases}\mathbf{A}_{00}=\left(-1\right)\mathbf{A}_{01}\mathbf{A}_{11}^{\top}\mathbf{A}_{10}.
\]
\[
\implies\left(\begin{array}{ccc}
\frac{\mathbf{A}_{00}}{\sqrt{2}} &  & \frac{\mathbf{A}_{01}}{\sqrt{2}}\\
\\
\frac{\mathbf{A}_{10}}{\sqrt{2}} &  & \frac{\mathbf{A}_{11}}{\sqrt{2}}
\end{array}\right)=\left(\begin{array}{ccc}
\frac{\left(-1\right)\mathbf{A}_{01}\mathbf{A}_{11}^{\top}\mathbf{A}_{10}}{\sqrt{2}} &  & \frac{\mathbf{A}_{01}}{\sqrt{2}}\\
\\
\frac{\mathbf{A}_{10}}{\sqrt{2}} &  & \frac{\mathbf{A}_{11}}{\sqrt{2}}
\end{array}\right).
\]
\textbf{Recall the canonical representation of the complex numbers
by $2\times2$} matrices described in Eq. (\ref{representation_of_CC}),
an unitary matrix $\mathbf{U}\in\mathbb{C}^{n\times n}$ can therefore
be seen as an $n\times n$ matrix of $2\times2$ block denoted $\mathbf{A}\in\left(\mathbb{R}^{2\times2}\right)^{n\times n}$
such that 
\[
\mathbf{A}\left[i,j\right]=\left(\begin{array}{ccc}
\Re\left(\mathbf{U}\left[i,j\right]\right) &  & -\Im\left(\mathbf{U}\left[i,j\right]\right)\\
\\
\Im\left(\mathbf{U}\left[i,j\right]\right) &  & \Re\left(\mathbf{U}\left[i,j\right]\right)
\end{array}\right).
\]
It follows that 
\[
\mathbf{U}\mathbf{U}^{*}=\mathbf{I}_{n}\Leftrightarrow\mathbf{A}\left(\mathbf{A}^{\top_{e}}\right)^{\top_{b}}=\mathbf{I}_{n}\otimes\mathbf{I}_{2}.
\]
It is therefore apparent that the algebra of complex numbers closely
relate to the algebra of $2\times2$ matrices and of real orthogonal
matrices in particular.

\subsection{The hypermatrix case.}

We now extend the discussion in section 7.1 to the hypermatrix case
to emphasize the compelling similarities. By analogy to the matrix
case, the action on the vector space $\mathbb{C}^{n\times1\times1}$
is specified in terms of a triple $\mathbf{A},\mathbf{B},\mathbf{C}\in\mathbb{C}^{n\times n\times n}$
as 
\[
\mathcal{T}_{\mathbf{A},\mathbf{B},\mathbf{C}}\,:\mathbb{C}^{n\times1\times1}\rightarrow\mathbb{C}^{n\times1\times1},\quad\mathbf{y}=\mathcal{T}_{\mathbf{A},\mathbf{B},\mathbf{C}}\left(\mathbf{x}\right),
\]
\[
\text{such that}
\]
\[
\forall\:0\le k<n,\quad\begin{cases}
\begin{array}{ccc}
\mathbf{y}\left[k\right] & = & \sqrt[3]{\text{Prod}_{\mathbf{P}_{k}}\left(\mathbf{x}^{\top^{2}},\mathbf{x}^{\top},\mathbf{x}\right)}.\\
\\
\mathbf{P}_{k} & = & \text{Prod}_{\boldsymbol{\Delta}^{(k)}}\left(\mathbf{A},\mathbf{B},\mathbf{C}\right)
\end{array} & \text{where }\end{cases}\boldsymbol{\Delta}^{(t)}\left[i,j,k\right]=\begin{cases}
\begin{array}{cc}
1 & \mbox{ if }\:0\le t=i=j<n\\
0 & \mbox{otherwise}
\end{array}\end{cases}.
\]
Invertibility in this context means that neither of the polynomials
in 
\[
0\ne\text{Resultant}_{\mathbf{x}}\left\{ \text{Prod}\left(\mathbf{x}^{\top^{2}},\mathbf{x}^{\top},\mathbf{x}\right):0\le k<n\right\} .
\]
is an identically non-zero constant. Recall that a triple $\mathbf{A},\mathbf{B},\mathbf{C}\in\mathbb{C}^{n\times n\times n}$
form an \emph{uncorrelated triple} if 
\[
\mbox{Prod}\left(\mathbf{A},\mathbf{B},\mathbf{C}\right)\left[i,j,k\right]=\begin{cases}
\begin{array}{cc}
1 & \text{if }i=j=k\\
0 & \text{otherwise}
\end{array}.\end{cases}.
\]
In the case where $\mathbf{A},\mathbf{B},\mathbf{C}\in\mathbb{C}^{n\times n\times n}$
form an uncorrelated triple, the map $\mathcal{T}_{\mathbf{A},\mathbf{B},\mathbf{C}}$
is subject to the resolution of identity 
\[
\mathbf{y}=\mathcal{T}_{\mathbf{A},\mathbf{B},\mathbf{C}}\left(\mathbf{x}\right)\implies\text{Prod}\left(\mathbf{y}^{\top^{2}},\mathbf{y}^{\top},\mathbf{y}\right)=\text{Prod}\left(\mathbf{x}^{\top^{2}},\mathbf{x}^{\top},\mathbf{x}\right).
\]
In other words, the map preserves the sum of cubes of the entries. 

In the case where $\mathbf{A},\mathbf{B},\mathbf{C}\in\mathbb{C}^{n\times n\times n}$
form an uncorrelated triple, $\mathbf{B}=\mathbf{A}^{\top^{2}}$ and
$\mathbf{C}=\mathbf{A}^{\top}$, the map $\mathcal{T}_{\mathbf{A},\mathbf{A}^{\top^{2}},\mathbf{A}^{\top}}$
is the third order hypermatrix analog of the vector isometry. This
latter observation therefore emphasizes the importance of orthogonal
hypermatrices. Recall that 
\[
\mathbf{X}\left[:,:,0\right]=\left(\begin{array}{rr}
x_{0} & x_{2}\\
x_{1} & x_{3}
\end{array}\right),\quad\mathbf{X}\left[:,:,1\right]=\left(\begin{array}{rr}
x_{4} & x_{6}\\
x_{5} & x_{7}
\end{array}\right),
\]
is orthogonal if
\[
\mbox{Prod}\left(\mathbf{X},\mathbf{X}^{\top^{2}},\mathbf{X}^{\top}\right)\left[:,:,0\right]=\left(\begin{array}{cc}
1 & 0\\
0 & 0
\end{array}\right),\quad\mbox{Prod}\left(\mathbf{X},\mathbf{X}^{\top^{2}},\mathbf{X}^{\top}\right)\left[:,:,1\right]=\left(\begin{array}{cc}
0 & 0\\
0 & 1
\end{array}\right).
\]
The corresponding constraints are therefore given by the polynomial
constraints
\[
\begin{cases}
\begin{array}{ccccc}
x_{1}x_{4}x_{5} & + & x_{3}x_{6}x_{7} & = & 0\\
x_{0}x_{1}x_{4} & + & x_{2}x_{3}x_{6} & = & 0\\
x_{0}^{3} & + & x_{2}^{3} & = & 1\\
x_{5}^{3} & + & x_{7}^{3} & = & 1
\end{array} & .\end{cases}
\]
When $0\ne\underset{0\le i<8}{\prod}x_{i}$, the above system of equations
yields the equivalence of
\[
\begin{array}{c}
\begin{cases}
\begin{array}{c}
0=x_{1}x_{4}x_{5}+x_{3}x_{6}x_{7}\\
\\
0=x_{0}x_{1}x_{4}+x_{2}x_{3}x_{6}
\end{array}\end{cases}\Leftrightarrow\begin{cases}
\begin{array}{c}
x_{1}x_{3}^{-1}x_{4}x_{5}x_{6}^{-1}x_{7}^{-1}=1\\
\\
x_{0}x_{1}x_{2}^{-1}x_{3}^{-2}x_{4}x_{6}^{-1}=1
\end{array}\end{cases}\end{array}
\]
\begin{equation}
\implies\left(\begin{array}{c}
x_{0}\\
x_{1}\\
x_{2}\\
x_{3}\\
x_{4}\\
x_{5}\\
x_{6}\\
x_{7}
\end{array}\right)=\left(\begin{array}{c}
\frac{v_{0}v_{3}}{v_{5}}\\
-\frac{v_{1}v_{4}v_{5}}{v_{2}v_{3}}\\
v_{0}\\
v_{1}\\
v_{2}\\
v_{3}\\
v_{4}\\
v_{5}
\end{array}\right)
\end{equation}

\[
\implies\mathbf{X}\left[:,:,0\right]=\left(\begin{array}{cc}
\frac{v_{0}v_{3}}{v_{5}} & v_{0}\\
-\frac{v_{1}v_{4}v_{5}}{v_{2}v_{3}} & v_{1}
\end{array}\right)\text{\; \ensuremath{\mathbf{X}\left[:,:,1\right]}=\ensuremath{\left(\begin{array}{cc}
v_{2} & v_{4}\\
v_{3} & v_{5}
\end{array}\right)}}.
\]
We account for the sum of cube constraints by normalizing appropriate
rows as follows

\[
\mathbf{X}\left[:,:,0\right]=\left(\begin{array}{ccc}
\frac{v_{0}v_{3}}{\sqrt[3]{v_{3}^{3}+v_{5}^{3}}} &  & \frac{v_{0}v_{5}}{\sqrt[3]{v_{3}^{3}+v_{5}^{3}}}\\
\\
-\frac{v_{1}v_{4}v_{5}}{v_{2}v_{3}} &  & v_{1}
\end{array}\right),\,\mathbf{X}\left[:,:,1\right]=\left(\begin{array}{ccc}
v_{2} &  & v_{4}\\
\\
\frac{v_{3}}{\sqrt[3]{v_{3}^{3}+v_{5}^{3}}} &  & \frac{v_{5}}{\sqrt[3]{v_{3}^{3}+v_{5}^{3}}}
\end{array}\right)
\]
where 
\[
v_{0}\in\left\{ \text{\ensuremath{\exp}\ensuremath{\left(\frac{2\pi\,k\sqrt{-1}}{3}\right)}}\,:\,0\le k<3\right\} 
\]

When $0=\underset{0\le i<8}{\prod}x_{i}$, The variables to be assigned
zero entries are indicated in the table below
\begin{center}
\begin{tabular}{|c|c|c|}
\hline 
$\left[x_{0}=0,x_{3}=0,x_{1}=0\right]$ & $\left[x_{0}=0,x_{6}=0,x_{1}=0\right]$ & $\left[x_{0}=0,x_{3}=0,x_{7}=0,x_{1}=0\right]$\tabularnewline
\hline 
$\left[x_{0}=0,x_{7}=0,x_{6}=0,x_{1}=0\right]$ & $\left[x_{0}=0,x_{4}=0,x_{3}=0\right]$ & $\left[x_{0}=0,x_{3}=0,x_{5}=0\right]$\tabularnewline
\hline 
$\left[x_{0}=0,x_{4}=0,x_{6}=0\right]$ & $\left[x_{0}=0,x_{4}=0,x_{3}=0,x_{7}=0\right]$ & $\left[x_{0}=0,x_{4}=0,x_{7}=0,x_{6}=0\right]$\tabularnewline
\hline 
$\left[x_{0}=0,x_{6}=0,x_{5}=0\right]$ & $\left[x_{3}=0,x_{2}=0,x_{1}=0\right]$ & $\left[x_{2}=0,x_{6}=0,x_{1}=0\right]$\tabularnewline
\hline 
$\left[x_{7}=0,x_{2}=0,x_{1}=0\right]$ & $\left[x_{3}=0,x_{1}=0\right]$ & $\left[x_{6}=0,x_{1}=0\right]$\tabularnewline
\hline 
$\left[x_{3}=0,x_{7}=0,x_{1}=0\right]$ & $\left[x_{7}=0,x_{6}=0,x_{1}=0\right]$ & $\left[x_{4}=0,x_{3}=0,x_{2}=0\right]$\tabularnewline
\hline 
$\left[x_{3}=0,x_{2}=0,x_{1}=0,x_{5}=0\right]$ & $\left[x_{4}=0,x_{3}=0,x_{2}=0,x_{5}=0\right]$ & $\left[x_{4}=0,x_{2}=0,x_{6}=0\right]$\tabularnewline
\hline 
$\left[x_{4}=0,x_{7}=0,x_{2}=0\right]$ & $\left[x_{2}=0,x_{6}=0,x_{1}=0,x_{5}=0\right]$ & $\left[x_{4}=0,x_{2}=0,x_{6}=0,x_{5}=0\right]$\tabularnewline
\hline 
$\left[x_{4}=0,x_{3}=0\right]$ & $\left[x_{3}=0,x_{1}=0,x_{5}=0\right]$ & $\left[x_{4}=0,x_{3}=0,x_{5}=0\right]$\tabularnewline
\hline 
$\left[x_{4}=0,x_{6}=0\right]$ & $\left[x_{4}=0,x_{3}=0,x_{7}=0\right]$ & $\left[x_{4}=0,x_{7}=0,x_{6}=0\right]$\tabularnewline
\hline 
$\left[x_{6}=0,x_{1}=0,x_{5}=0\right]$ & $\left[x_{4}=0,x_{6}=0,x_{5}=0\right]$ & \tabularnewline
\hline 
\end{tabular}
\par\end{center}

To express some important invariants of orthogonal hypermatrices,
we extend the index rotation operation to third order hypermatrices
and is denoted by $\mathbf{A}^{\text{R}_{\left[\theta_{x},\theta_{y},\theta_{z}\right]}}$
for $\theta_{x},\theta_{y},\theta_{z}\in\left\{ 0\cdot\frac{2\pi}{4},1\cdot\frac{2\pi}{4},2\cdot\frac{2\pi}{4},3\cdot\frac{2\pi}{4}\right\} $
such that $\mathbf{A}^{\text{R}_{\left[\theta_{x}0,0\right]}}$ denotes
the hypermatrix which result from performing the index rotation by
angle $\theta_{x}$ to each row slices of $\mathbf{A}$. Similarly,
$\mathbf{A}^{\text{R}_{\left[0,\theta_{y},0\right]}}$ denotes the
hypermatrix which result from performing the index rotation by angle
$\theta_{y}$ to each column slice of $\mathbf{A}$ and finally $\mathbf{A}^{\text{R}_{\left[0,0,\theta_{z}\right]}}$
denotes the hypermatrix which result from performing the index rotation
by angle $\theta_{z}$ to each depth slice of $\mathbf{A}$. The index
rotation $\mathbf{A}^{\text{R}_{\left[\theta_{x},\theta_{y},\theta_{z}\right]}}$
is performed relative to the axis $x$, $y$ and $z$ in that order.
For instance we have
\[
\left(\mathbf{A}^{\text{R}_{\left[0,0,\theta_{z}\right]}}\right)\left[i,j,k\right]=
\]
\[
\mathbf{A}\left[\left(i-\frac{n-1}{2}\right)\cos\theta+\left(\frac{n-1}{2}-j\right)\sin\theta+\frac{n-1}{2},\,\left(i-\frac{n-1}{2}\right)\sin\theta-\left(\frac{n-1}{2}-j\right)\cos\theta+\frac{n-1}{2},\,k\right].
\]
\begin{equation}
\mbox{Prod}\left(\mathbf{X},\mathbf{X}^{\top^{2}},\mathbf{X}^{\top}\right)=\Delta_{2}\implies\mbox{Prod}\left(\mathbf{X}^{\text{R}_{\left[\frac{2\pi k_{0}}{4},\frac{2\pi k_{1}}{4},\frac{2\pi k_{2}}{4}\right]}},\left(\mathbf{X}^{\text{R}_{\left[\frac{2\pi k_{0}}{4},\frac{2\pi k_{1}}{4},\frac{2\pi k_{2}}{4}\right]}}\right)^{\top^{2}},\left(\mathbf{X}^{\text{R}_{\left[\frac{2\pi k_{0}}{4},\frac{2\pi k_{1}}{4},\frac{2\pi k_{2}}{4}\right]}}\right)^{\top}\right)=\boldsymbol{\Delta}
\end{equation}
where $\left[k_{0}\frac{2\pi}{4},k_{1}\frac{2\pi}{4},k_{2}\frac{2\pi}{4}\right]$
belong to values indicated in the table below 
\begin{center}
\begin{tabular}{|c|c|c|c|}
\hline 
$\left[0,0,0\right]$ & $\left[0,0,\frac{3}{2}\,\pi\right]$ & $\left[0,\frac{1}{2}\,\pi,0\right]$ & $\left[0,\frac{1}{2}\,\pi,\pi\right]$\tabularnewline
\hline 
$\left[0,\pi,\frac{1}{2}\,\pi\right]$ & $\left[0,\pi,\pi\right]$ & $\left[0,\frac{3}{2}\,\pi,\frac{1}{2}\,\pi\right]$ & $\left[0,\frac{3}{2}\,\pi,\frac{3}{2}\,\pi\right]$\tabularnewline
\hline 
$\left[\frac{1}{2}\,\pi,0,\frac{1}{2}\,\pi\right]$ & $\left[\frac{1}{2}\,\pi,0,\pi\right]$ & $\left[\frac{1}{2}\,\pi,\frac{1}{2}\,\pi,\frac{1}{2}\,\pi\right]$ & $\left[\frac{1}{2}\,\pi,\frac{1}{2}\,\pi,\frac{3}{2}\,\pi\right]$\tabularnewline
\hline 
$\left[\frac{1}{2}\,\pi,\pi,0\right]$ & $\left[\frac{1}{2}\,\pi,\pi,\frac{3}{2}\,\pi\right]$ & $\left[\frac{1}{2}\,\pi,\frac{3}{2}\,\pi,0\right]$ & $\left[\frac{1}{2}\,\pi,\frac{3}{2}\,\pi,\pi\right]$\tabularnewline
\hline 
$\left[\pi,0,0\right]$ & $\left[\pi,0,\frac{3}{2}\,\pi\right]$ & $\left[\pi,\frac{1}{2}\,\pi,0\right]$ & $\left[\pi,\frac{1}{2}\,\pi,\pi\right]$\tabularnewline
\hline 
$\left[\pi,\pi,\frac{1}{2}\,\pi\right]$ & $\left[\pi,\pi,\pi\right]$ & $\left[\pi,\frac{3}{2}\,\pi,\frac{1}{2}\,\pi\right]$ & $\left[\pi,\frac{3}{2}\,\pi,\frac{3}{2}\,\pi\right]$\tabularnewline
\hline 
$\left[\frac{3}{2}\,\pi,0,\frac{1}{2}\,\pi\right]$ & $\left[\frac{3}{2}\,\pi,0,\pi\right]$ & $\left[\frac{3}{2}\,\pi,\frac{1}{2}\,\pi,\frac{1}{2}\,\pi\right]$ & $\left[\frac{3}{2}\,\pi,\frac{1}{2}\,\pi,\frac{3}{2}\,\pi\right]$\tabularnewline
\hline 
$\left[\frac{3}{2}\,\pi,\pi,0\right]$ & $\left[\frac{3}{2}\,\pi,\pi,\frac{3}{2}\,\pi\right]$ & $\left[\frac{3}{2}\,\pi,\frac{3}{2}\,\pi,0\right]$ & $\left[\frac{3}{2}\,\pi,\frac{3}{2}\,\pi,\pi\right]$\tabularnewline
\hline 
\end{tabular}
\par\end{center}

As shown in \cite{Gnang2017238} if Prod$\left(\mathbf{X},\mathbf{X}^{\top^{2}},\mathbf{X}^{\top}\right)=\boldsymbol{\Delta}=$
Prod$\left(\mathbf{Y},\mathbf{Y}^{\top^{2}},\mathbf{Y}{}^{\top}\right)$
then we have 

\[
\text{Prod}\left(\left(\mathbf{X}\oplus\mathbf{Y}\right),\left(\mathbf{X}\oplus\mathbf{Y}\right)^{\top^{2}},\left(\mathbf{X}\oplus\mathbf{Y}\right)^{\top}\right)=\boldsymbol{\Delta}\oplus\boldsymbol{\Delta}
\]
\[
\text{and}
\]
\[
\text{Prod}\left(\left(\mathbf{X}\otimes\mathbf{Y}\right),\left(\mathbf{X}\otimes\mathbf{Y}\right)^{\top^{2}},\left(\mathbf{X}\otimes\mathbf{Y}\right)^{\top}\right)=\boldsymbol{\Delta}\otimes\boldsymbol{\Delta}
\]

Consider block operation of hymatrices 
\begin{equation}
\mathcal{A}\left[:,:,0\right]=\left(\begin{array}{ccc}
\mathbf{A}_{000} &  & \mathbf{A}_{010}\\
\\
\mathbf{A}_{100} &  & \mathbf{A}_{110}
\end{array}\right)\mathcal{A}\left[:,:,1\right]=\left(\begin{array}{ccc}
\mathbf{A}_{001} &  & \mathbf{A}_{011}\\
\\
\mathbf{A}_{101} &  & \mathbf{A}_{111}
\end{array}\right)
\end{equation}

\begin{equation}
\mathcal{A}{}^{\top_{b}^{t}}=\sum_{0\le i,j,k<2}\text{Prod}\left(\mathbf{K}_{0}\left[:,i,:\right],\mathbf{K}_{1}\left[:,:,j\right],\mathbf{K}_{2}\left[k,:,:\right]\right)^{\top^{t}}\otimes\mathbf{A}_{ijk}
\end{equation}

\[
\text{s.t.}
\]

\[
\mathbf{K}_{0}[:,:,0]=\left(\begin{array}{rr}
1 & 0\\
0 & 1
\end{array}\right);\:\mathbf{K}_{0}[:,:,1]=\left(\begin{array}{rr}
1 & 0\\
0 & 1
\end{array}\right)
\]

\[
\mathbf{K}_{1}[:,:,0]=\left(\begin{array}{rr}
1 & 0\\
1 & 0
\end{array}\right);\:\mathbf{K}_{1}[:,:,1]=\left(\begin{array}{rr}
0 & 1\\
0 & 1
\end{array}\right)
\]

\[
\mathbf{K}_{2}[:,:,0]=\left(\begin{array}{rr}
1 & 1\\
0 & 0
\end{array}\right);\:\mathbf{K}_{2}[:,:,1]=\left(\begin{array}{rr}
0 & 0\\
1 & 1
\end{array}\right)
\]
\begin{equation}
\mathbf{A}{}^{\top_{e}^{t}}=\sum_{0\le i,j,k<2}\text{Prod}\left(\mathbf{K}_{0}\left[:,i,:\right],\mathbf{K}_{1}\left[:,:,j\right],\mathbf{K}_{2}\left[k,:,:\right]\right)\otimes\mathbf{A}_{ijk}^{\top^{t}}
\end{equation}
\begin{equation}
\mathbf{A}^{\text{R}_{_{b},\frac{2\pi k}{4}}}=\sum_{0\le i,j,k<2}\text{Prod}\left(\mathbf{K}_{0}\left[:,i,:\right],\mathbf{K}_{1}\left[:,:,j\right],\mathbf{K}_{2}\left[k,:,:\right]\right)^{\text{R}_{\frac{2\pi k}{4}}}\otimes\mathbf{A}_{ijk}
\end{equation}
\begin{equation}
\mathbf{A}^{\text{R}_{e,\frac{2\pi k}{4}}}=\sum_{0\le i,j,k<2}\text{Prod}\left(\mathbf{K}_{0}\left[:,i,:\right],\mathbf{K}_{1}\left[:,:,j\right],\mathbf{K}_{2}\left[k,:,:\right]\right)\otimes\mathbf{A}_{ijk}^{\text{R}_{\frac{2\pi k}{4}}}
\end{equation}
Similarly to the matrix case, if $\mathbf{A}$ is block hypermatrix
whose invidividual blocks are orthogonal hypermatrices all of the
same size and all subject to 
\[
\text{Prod}\left(\mathbf{X},\mathbf{X}^{\top^{2}},\mathbf{X}^{\top}\right)=\boldsymbol{\Delta}
\]
 then it follows that 
\[
\forall\:0\le i<n,\quad\text{Prod}\left(\frac{\mathbf{A}}{\sqrt[3]{n}},\left(\frac{\mathbf{A}}{\sqrt[3]{n}}^{\top_{e}^{2}}\right)^{\top_{b}^{2}},\left(\frac{\mathbf{A}}{\sqrt[3]{n}}^{\top_{e}}\right)^{\top_{b}}\right)\left[i,i,i\right]=\boldsymbol{\Delta}
\]
In which case 
\[
\text{Prod}\left(\frac{\mathbf{A}}{\sqrt[3]{n}},\left(\frac{\mathbf{A}}{\sqrt[3]{n}}^{\top_{e}^{2}}\right)^{\top_{b}^{2}},\left(\frac{\mathbf{A}}{\sqrt[3]{n}}^{\top_{e}}\right)^{\top_{b}}\right)=\mathbf{0}\implies\mathbf{A}=\mathbf{0}.
\]
In the case of $2\times2\times2$ block hypermatrix 
\[
\mathbf{A}[:,:,0]=\frac{1}{\sqrt[3]{2}}\left(\begin{array}{rr}
\mathbf{A}_{000} & \mathbf{A}_{010}\\
\mathbf{A}_{100} & \mathbf{A}_{110}
\end{array}\right),\quad\mathbf{A}[:,:,1]=\frac{1}{\sqrt[3]{2}}\left(\begin{array}{rr}
\mathbf{A}_{001} & \mathbf{A}_{011}\\
\mathbf{A}_{101} & \mathbf{A}_{111}
\end{array}\right)
\]
where 
\[
\forall\,\mathbf{X}\in\left\{ \mathbf{A}_{000},\mathbf{A}_{100},\mathbf{A}_{010},\mathbf{A}_{110},\mathbf{A}_{001},\mathbf{A}_{101},\mathbf{A}_{011},\mathbf{A}_{111}\right\} ,
\]
we have 
\[
\text{Prod}\left(\mathbf{X},\mathbf{X}^{\top^{2}},\mathbf{X}^{\top}\right)=\boldsymbol{\Delta}
\]
is expressed by
\[
\text{Prod}\left(\mathbf{A},\left(\mathbf{A}^{\top_{e}^{2}}\right)^{\top_{b}^{2}},\left(\mathbf{A}^{\top_{e}}\right)^{\top_{b}}\right)[:,:,0]=
\]
 
\[
\left(\begin{array}{ccc}
\boldsymbol{\Delta} &  & \frac{\text{Prod}\left(\mathbf{A}_{000},\mathbf{A}_{100}^{\top^{2}},\mathbf{A}_{001}^{\top}\right)+\text{Prod}\left(A_{010},\mathbf{A}_{110}^{\top^{2}},\mathbf{A}_{011}^{\top}\right)}{2}\\
\\
\frac{\text{Prod}\left(\mathbf{A}_{100},\mathbf{A}_{001}^{\top^{2}},\mathbf{A}_{000}^{\top}\right)+\text{Prod}\left(\mathbf{A}_{110},\mathbf{A}_{011}^{\top^{2}},\mathbf{A}_{010}^{\top}\right)}{2} &  & \frac{\text{Prod}\left(\mathbf{A}_{100},\mathbf{A}_{101}^{\top^{2}},\mathbf{A}_{001}^{\top}\right)+\text{Prod}\left(\mathbf{A}_{110},\mathbf{A}_{111}^{\top^{2}},\mathbf{A}_{011}^{\top}\right)}{2}
\end{array}\right),
\]
\[
\text{Prod}\left(\mathbf{A},\left(\mathbf{A}^{\top_{e}^{2}}\right)^{\top_{b}^{2}},\left(\mathbf{A}^{\top_{e}}\right)^{\top_{b}}\right)[:,:,1]=
\]
\[
\left(\begin{array}{ccc}
\frac{\text{Prod}\left(\mathbf{A}_{001},\mathbf{A}_{000}^{\top^{2}},\mathbf{A}_{100}^{\top}\right)+\text{Prod}\left(A_{011},\mathbf{A}_{010}^{\top^{2}},\mathbf{A}_{110}^{\top}\right)}{2} &  & \frac{\text{Prod}\left(\mathbf{A}_{001},\mathbf{A}_{100}^{\top^{2}},\mathbf{A}_{101}^{\top}\right)+\text{Prod}\left(\mathbf{A}_{011},\mathbf{A}_{110}^{\top^{2}},\mathbf{A}_{111}^{\top}\right)}{2}\\
\\
\frac{\text{Prod}\left(A_{101},\mathbf{A}_{001}^{\top^{2}},\mathbf{A}_{100}^{\top}\right)+\text{Prod}\left(\mathbf{A}_{111},\mathbf{A}_{011}^{\top^{2}},\mathbf{A}_{110}^{\top}\right)}{2} &  & \boldsymbol{\Delta}
\end{array}\right).
\]
A necessary condition for the resulting block hypermatrix to be orthogonal
is specified by the constraints 
\[
\begin{cases}
\begin{array}{ccc}
\text{Prod}\left(\mathbf{A}_{000},\mathbf{A}_{100}^{\top^{2}},\mathbf{A}_{001}^{\top}\right)+\text{Prod}\left(A_{010},\mathbf{A}_{110}^{\top^{2}},\mathbf{A}_{011}^{\top}\right) & = & \mathbf{0}\\
\\
\text{Prod}\left(\mathbf{A}_{100},\mathbf{A}_{001}^{\top^{2}},\mathbf{A}_{000}^{\top}\right)+\text{Prod}\left(\mathbf{A}_{110},\mathbf{A}_{011}^{\top^{2}},\mathbf{A}_{010}^{\top}\right) & = & \mathbf{0}\\
\\
\text{Prod}\left(\mathbf{A}_{100},\mathbf{A}_{101}^{\top^{2}},\mathbf{A}_{001}^{\top}\right)+\text{Prod}\left(\mathbf{A}_{110},\mathbf{A}_{111}^{\top^{2}},\mathbf{A}_{011}^{\top}\right) & = & \mathbf{0}\\
\\
\text{Prod}\left(\mathbf{A}_{001},\mathbf{A}_{000}^{\top^{2}},\mathbf{A}_{100}^{\top}\right)+\text{Prod}\left(A_{011},\mathbf{A}_{010}^{\top^{2}},\mathbf{A}_{110}^{\top}\right) & = & \mathbf{0}\\
\\
\text{Prod}\left(\mathbf{A}_{001},\mathbf{A}_{100}^{\top^{2}},\mathbf{A}_{101}^{\top}\right)+\text{Prod}\left(\mathbf{A}_{011},\mathbf{A}_{110}^{\top^{2}},\mathbf{A}_{111}^{\top}\right) & = & \mathbf{0}\\
\\
\text{Prod}\left(A_{101},\mathbf{A}_{001}^{\top^{2}},\mathbf{A}_{100}^{\top}\right)+\text{Prod}\left(\mathbf{A}_{111},\mathbf{A}_{011}^{\top^{2}},\mathbf{A}_{110}^{\top}\right) & = & \mathbf{0}
\end{array}\end{cases}
\]
\providecommand{\bysame}{\leavevmode\hbox to3em{\hrulefill}\thinspace}
\providecommand{\MR}{\relax\ifhmode\unskip\space\fi MR }
\providecommand{\MRhref}[2]{%
  \href{http://www.ams.org/mathscinet-getitem?mr=#1}{#2}
}
\providecommand{\href}[2]{#2}

\end{document}